\documentclass{article}

\usepackage{amsmath}
\usepackage{amsfonts}
\usepackage{amssymb}
\usepackage{amsthm}

\usepackage[greek,french,english]{babel}
\usepackage[T1]{fontenc}
\usepackage[latin1]{inputenc}
\usepackage{graphics}
\usepackage{vmargin}
\usepackage{enumerate}

\setmarginsrb{3cm}{3cm}{3cm}{3cm}{0cm}{0cm}{0cm}{1.5cm}

\newcommand{\R}{\mathbb R}
\newcommand{\N}{\mathbb N}
\newcommand{\Z}{\mathbb Z}
\newcommand{\T}{\mathbb T}
\newcommand{\e}{\varepsilon}
\newcommand{\be}{\begin{equation}}
\newcommand{\ee}{\end{equation}}
\newcommand{\dv}{\mathrm{div}}
\newcommand{\p}{\partial}

\newcommand{\tU}{\tilde{U}}
\newcommand{\tW}{\tilde{W}}
\newcommand{\ds}{\displaystyle}

\newcommand{\mean}[1]{\left\langle #1\right\rangle}

\newcommand{\cA}{\mathcal A}

\newcommand{\Uapp}{U^\text{app}}
\newcommand{\bV}{\bar V}
\newcommand{\unde}{U_n^{\delta,\e}}

\newtheorem{prop}{Proposition}[section]
\newtheorem{lem}{Lemma}[section]

\newtheorem{thm}{Theorem}
\newtheorem{defi}{Definition}[section]
\newtheorem{corol}{Corollary}[section]

\theoremstyle{definition}
\newtheorem{rem}{Remark}

\begin{document}
\title{Stability of periodic stationary solutions of scalar conservation laws with space-periodic flux}

\author{Anne-Laure Dalibard\thanks{CNRS, UMR 8553-DMA, Ecole Normale Sup\'erieure, 45 rue d'Ulm, 75005 Paris, FRANCE. e-mail: {\tt Anne-Laure.Dalibard@ens.fr}}
}

\bibliographystyle{amsplain}
\maketitle

\begin{abstract}
This article investigates the long-time behaviour of parabolic scalar conservation laws of the type $\partial_t u + \mathrm{div}_yA(y,u) - \Delta_y u=0$, where $y\in\mathbb R^N$ and the flux $A$ is periodic in $y$. More specifically, we consider the case when the initial data is an $L^1$ disturbance of a stationary periodic solution. We show, under polynomial growth assumptions on the flux, that the difference between $u$ and the stationary solution vanishes for large times in $L^1$ norm. The proof uses a self-similar change of variables which is well-suited for the analysis of the long time behaviour of parabolic equations. Then, convergence in self-similar variables follows from arguments from dynamical systems theory. One crucial point is to obtain  compactness in $L^1$ on the family of rescaled solutions; this is achieved by deriving uniform bounds in weighted $L^2$ spaces.
\end{abstract}



\section{Introduction}
The goal of this article is to study the long time limit of solutions of the equation
\be
\p_t u + \dv_y A(y,u) -\Delta_y u=0,\quad t>0,\ y\in\R^N,\label{LCS}
\ee
where the flux $A:\R^N\times \R\to \R^N$ is assumed to be $\T^N$-periodic with respect to its first variable. Here and in the rest of the article, $\T^N$ denotes the $N$-dimensional torus, i.e. $\T^N=(\R/\Z)^N$.

Classical results on scalar conservation laws (see for instance \cite{SerreLCS,kruzhkov}) ensure that the semi-group associated with equation \eqref{LCS} is well-defined in $L^1(\R^N) + L^\infty(\R^N)$. The case when the initial data belongs to $U(y) + L^1(\R)$ (when $N=1$), where $U$ is a viscous shock profile of equation \eqref{LCS} has already been dealt with in a previous article, see \cite{longtime}. In the present paper, we restrict our study to the case when the inital data belongs to $v(y) + L^1(\R^N)$, where  $v$ is a given periodic stationary solution of \eqref{LCS}.

When the flux $A$ is linear, say
$$
A(y,u)=\alpha(y) u,
$$
this study coincides (at least for some particular functions $\alpha$) with the one led by Adrien Blanchet, Jean Dolbeault and Michal Kowalczyk in \cite{BDK} on the large time behaviour of Brownian ratchets, as we will explain in Remark \ref{rmk:ratchets}. It is proved in \cite{BDK}  that if the flux $A$ is linear and if
\be
\limsup_{t\to\infty}\frac{1}{(1+2t)^2}\int_{\R^N}|u(t,y) - v(y)|\; (y-c t)^4\:dy<\infty\label{hyp:moments}
\ee
for some velocity $c\in\R^N$ which will be defined later on (see \eqref{def:c}), then there exists a constant $C$ and a number $\kappa\in(0,1/2)$ such that
\be
\int_{\R^N}\left| u(t,y) -v(y) - \frac{M f_0(y)}{(1+2t)^{N/2}}F \left( \frac{y-ct}{\sqrt{1+2t}} \right)  \right|dy\leq C t^{-\kappa},\label{in:BDK}
\ee
where $f_0$ is the  solution of an elliptic equation in $\T^N$ (see \eqref{eq:f0}), $F$ is a Gaussian profile, and $M$ is the mass of the initial disturbance, i.e.
$$
M=\int_{\R^N} (u_{|t=0}- v).
$$
Unfortunately, as we explained in \cite{longtime}, the above result does not imply that the same convergence holds in the nonlinear case. Moreover, the proof of \cite{BDK}, which is based on entropy dissipation methods together with Log-Sobolev Poincar\'e inequalities, can hardly be transposed as such to a nonlinear setting, although attempts in this direction have been made, see for instance \cite{DFM}.
Hence we have chosen here a slightly different approach, which enables us to extend some of the results of \cite{BDK} to a nonlinear context. Additionally, we recover a weaker version of the convergence \eqref{in:BDK}, but without the need for assumption \eqref{hyp:moments}. In fact, we prove that \eqref{hyp:moments} holds for a large class of initial data in $v+ L^1(\R^N)$.

The present work is also embedded in the broader study of the long time behaviour of conservation laws. We refer the interested reader to the review paper by D. Serre \cite{SerreHandbook} (and the references therein) for a thorough description of the homogeneous case, in which the author investigates  the stability of stationary solutions of scalar conservation laws in various models (parabolic and hyperbolic settings, relaxation models...)

Before stating the main results of this paper, let us now recall a few properties of equation \eqref{LCS}. First, according to a result of \cite{homogpara}, periodic stationary solutions of \eqref{LCS} exist, provided the flux $A$ satisfies some growth assumptions. In fact, several different growth regimes were studied in \cite{homogpara}; we only recall one of them here, which is the most relevant with regards to our purposes. In the rest of the article, we assume that $A$ belongs to $W^{1,\infty}_\text{loc}(\T^N\times \R)^N$, and that 
\be
\exists p_0\in \R,\ \forall y\in\T^N,\ \dv_y A(y,p_0)=
0.\label{hyp:A1}
\ee
We also assume that
there exists $n\in(0, (N+2)/N)$ such that
\begin{multline}
\label{hyp:A2}
\forall P>0, \ \exists C_P>0,\ \forall (p,q)\in\R^2,\ |p|\leq P,\\ \left\{ 
\begin{array}{l}
\left| \p_p A(y,p+q)- \p_p A(y,p)  \right|\leq C_P (|q| + |q|^{n}),\\
\left| \dv_y A(y,p+q)- \dv_y A(y,p)  \right|\leq C_P (|q| + |q|^{n}).
\end{array}
\right.	
\end{multline}
These assumptions were introduced in \cite{longtime,homogpara}. They ensure that for any $q\in\R$, there exists a unique periodic stationary solution of \eqref{LCS} with mean value $q$; we refer to \cite{homogpara} for a discussion of the optimality of conditions \eqref{hyp:A1}, \eqref{hyp:A2}.  Moreover, if $u$ is a solution of \eqref{LCS} with initial data
 $u_{|t=0}\in v + L^1\cap L^\infty(\R^N)$, where $v\in W^{1,\infty}(\R^N)$ is any stationary solution of \eqref{LCS}, then $u\in L^\infty([0,\infty)\times\R^N)$. This result will be used several times in the article, and its proof is recalled in Appendix A.

\vskip1mm
We now introduce the profiles which characterize the asymptotic behaviour of the function $u$; first, the function $f_0$ occurring in \eqref{in:BDK} is the unique solution in $H^1(\T^N)$ of the equation
\be\label{eq:f0}
-\Delta_y f_0 + \dv_y ( \alpha_1 f_0)=0,\quad  \mean{f_0}=1,
\ee
where $\alpha_1(y):=(\p_p A) (y,v(y))\in L^\infty(\T^N)^N.$ Above and in the rest of the article, the notation $\mean{\cdot}$ stands for the average on the torus $\T^N$, that is
$$
\mean{f}:=\int_{\T^N} f.
$$
The drift velocity $c$ is then defined by 
\be\label{def:c}
c=N\mean{\alpha_1 f_0}.
\ee
The last function which will appear in the asymptotic profile of $u$ is the equivalent, in the non-linear case, of the Gaussian profile $F$ occurring in \eqref{in:BDK}; it is the unique solution, in a suitable functional space, of an elliptic equation of the form
$$
-\sum_{1\leq i,j\leq N}\eta_{i,j} \p_i\p_j F_M - \dv_x(x F_M) + a\cdot \nabla_x F_M^2=0\text{ in } \R^N,\text{ with }\int_{\R^N} F_M=M\in\R,
$$
where the coefficients $\eta_{i,j}$ and $a$ are constant, and the matrix $(\eta_{i,j})_{1\leq i,j\leq N}$ is coercive. Unfortunately, giving the precise definition of $\eta_{i,j}$ and $a$ would take us too far at this stage. We merely recall that thanks to a result of J. Aguirre, M. Escobedo, and E. Zuazua (see \cite{AEZ}), the above equation has a unique solution for all $M\in\R$, and we refer to the next section for more details.

\vskip1mm

The main result of this paper is the following:
\begin{thm}
Let $A\in W^{5,\infty}_\text{loc}(\T^N\times \R)^N$, and assume that $A$ satisfies \eqref{hyp:A1}, \eqref{hyp:A2}.

Let $v$ be a periodic stationary solution of \eqref{LCS}, and let $u_{ini}\in v + L^1(\R^N)$. Let $u$ be the unique solution of \eqref{LCS} with initial condition $u_{|t=0}=u_{ini}.$
Set 
$$
M:=\int_{\R} (u_{ini}-v)\:dy.
$$
Then as $t\to\infty$,
$$
\int_{\R^N} \left|u(t,y) - v(y) - \frac{1}{(1+2t)^{N/2}}f_0(y) F_M \left(\frac{y-c t}{\sqrt{1+2t}}\right)   \right| dy\to 0.
$$

\label{thm:conv}
\end{thm}

\begin{rem}
In fact, the regularity assumptions on the flux $A$ are not as stringent as stated in the Theorem above. In particular, the conditions on the derivatives with respect to the space variable $y$ can be considerably reduced. When looking closely at the proof, the correct regularity assumptions on $A$ are
$$
\begin{aligned}
\p_p^k A\in L^\infty_\text{loc}(\T^N\times \R)^N\quad\forall k\in\{0,1,\cdots, 4 \},\\
\dv_y A,\dv_y \p_p^4 A \in L^\infty_\text{loc}(\T^N\times \R).
\end{aligned}
$$
\end{rem}

\begin{rem}
Let us now make precise the link between brownian ratchets and equation \eqref{LCS} in the linear case. In \cite{BDK}, A. Blanchet, J. Dolbeault and M. Kowalczyk study the long time behaviour of the solution $f=f(t,y)$ of the equation
\be
\p_t f = \Delta_y f + \dv_y(\nabla \psi(y-\omega t) f),\quad t>0,\ y\in\R^N,\label{eq:BDK}
\ee
with $\psi\in\mathcal C^2(\T^N),$ $\omega\in\R^N$. Setting
$$
u(t,y)=f(t,y+\omega t)\quad \forall t>0,\ \forall y\in\R^N,
$$
we see that $u$ satisfies
$$
\p_t u + \dv_y (\alpha(y) u ) - \Delta_y u=0,
$$
where the drift coefficient $\alpha$ is given by 
\be\label{alpha_ratchet}
\alpha(y)= - \omega - \nabla_y \psi(y).
\ee
Hence the study of \eqref{eq:BDK} and that of \eqref{LCS} in the linear case are closely related; they are strictly equivalent in dimension one, since any function $\alpha\in\mathcal C^1(\T)$ can be decomposed as 
$$
\alpha= \int_\T \alpha + \left(\alpha - \int_\T \alpha\right)=  \int_\T \alpha + \p_y \phi, \quad\text{for some }\phi\in\mathcal C^2(\T).
$$
The equivalence does not hold when $N\geq 2$, but in fact, all the results of \cite{BDK} remain true for an arbitrary drift $\alpha\in\mathcal C^1(\T^N)$ (using exactly the same techniques as the ones developed in \cite{BDK}). This will be a consequence of the analysis  we will perform in the next sections. The choice for a  function $\alpha$ with the structure \eqref{alpha_ratchet} stems from physical considerations (see \cite{BDK2}): equation \eqref{eq:BDK} describes the evolution of the density of particles in a traveling potential, moving with constant speed $\omega$.

\label{rmk:ratchets}
\end{rem}

In the course of the proof of Theorem \ref{thm:conv}, we will also prove that condition \eqref{hyp:moments} holds for a large class of initial data. The precise result is the following:

\begin{prop}\label{prop:linear}
Assume that the flux $A$ is linear, and that $u_{ini}\in v+ L^1(\R^N)$ is such that
$$
\exists m>N+8,\quad \int_{\R^N} |u_{ini}(y) - v(y)|^2 (1+ |y|^2)^{m/2}\:dy<\infty.
$$ 
Let $u$ be the unique solution of \eqref{LCS} with initial data $u_{ini}$. Then \eqref{hyp:moments} is satisfied. As a consequence (see \cite{BDK}), \eqref{in:BDK} holds.
\end{prop}
Hence for linear fluxes and for a large range of initial data, a rate of convergence can be given. The derivation of convergence rates in the non-linear case goes beyond the scope of this article; in fact, the standard methods to derive convergence rates rely on the use of entropy-entropy dissipation inequalities (see \cite{DFM} in the case of the Burgers equation), which we have chosen not to use here.

\vskip1mm

Another consequence of Theorem \ref{thm:conv} is the stability of stationary shock profiles of equation \eqref{LCS} (see \cite{longtime}) in dimension one: a stationary shock profile is a stationary solution of \eqref{LCS} with $N=1$, which is asymptotic as $y\to\pm \infty$ to periodic stationary solutions of \eqref{LCS}. It was proved in \cite{longtime} that the stability of shock profiles is a consequence of the stability of periodic stationary solutions. Thus we have the following 
\begin{corol}
Assume that $N=1$, and that the hypotheses of Theorem \ref{thm:conv} are satisfied. Let $U\in L^\infty(\R)$ be a stationary shock profile of \eqref{LCS}. Let $u_{ini}\in U+ L^1(\R)$ such that
$$
\int_{\R}(u_{ini} - U)=0,
$$
and let $u$ be the unique solution of \eqref{LCS} with initial data $u_{ini}$. Then
$$
\lim_{t\to\infty}\| u(t)- U\|_{L^1(\R)}=0.
$$
\end{corol}

\vskip1mm

The strategy of proof of Theorem \ref{thm:conv} is close to the one developed in \cite{EZ}, in which M. Escobedo and E. Zuazua study the long time behaviour of a homogeneous version of \eqref{LCS}; we also refer the interested reader to \cite{EVZ}, in which M. Escobedo, J.L. Vazquez and E. Zuazua extend the analysis performed in \cite{EZ} to the case when the flux has sub-critical growth. The first step of the analysis consists in a self-similar change of variables, which helps us to focus on the appropriate length scales; this will be done in the next section, in which we also derive the equations on the limit profiles $f_0$ and $F_M$. Then, in section \ref{sec:bounds}, we obtain some compactness on the rescaled sequence by deriving some uniform $L^2$ bounds in weighted spaces. Eventually, we conclude the proof in Section~\ref{sec:conclusion} by using semi-group arguments inherited from dynamical systems theory.

\vskip1mm
Throughout the article, we will use the following notation: if $\psi\in L^\infty_\text{loc}(\R^N)$, we set, for all $p\in[1,\infty)$,
$$
\begin{aligned}
L^p (\psi):=\left\{ u\in L^p_\text{loc}(\R^N),\ \int_{\R^N} |u|^p\psi <+\infty\right\},\\
\text{and }\|u\|_{L^p(\psi)}=\left( \int_{\R^N} |u|^p\psi \right)^{1/p},\\
H^1(\psi):=\left\{ u\in L^2(\psi),\ \nabla u \in L^2(\psi)\right\},\\
\text{and }\|u\|_{H^1(\psi)}^2= \|u\|_{L^2(\psi)}^2 + \|\nabla u\|_{L^2(\psi)}^2.
\end{aligned}
$$
Sobolev spaces of the type $W^{s,p}(\psi)$, $H^s(\psi)$, with $s\in\N$ arbitrary and $p\in[1,\infty)$, are defined in a similar fashion. When we write $\|u\|_p$, or $\|u\|_{L^p}$, without specifying a weight function, we always refer to the usual $L^p$ norm in $\R^N$, with respect to the Lebesgue measure (i.e. $\psi\equiv 1$).

\section{The homogenized system}

The goal of this section is to analyze the expected asymptotic behaviour of the solution $u(t)$ of equation \eqref{LCS}; to that end, we change the space and time variables and introduce a parabolic scaling, which is appropriate for  the study of the long time behaviour of diffusion equations. Then, using a two-scale Ansatz in space and time which was introduced in \cite{BDK}, we construct an approximate solution of the rescaled system. Eventually, we recall and derive several properties of the limit system.

\subsection{Parabolic scaling}

Consider the solution $u\in L^\infty_\text{loc}([0,\infty)\times \R^N)$ of \eqref{LCS}, with $u_{|t=0}=u_{ini}\in v + L^1\cap L^\infty(\R^N)$. It is a classical feature of scalar conservation laws  that the semi-group associated with \eqref{LCS} is contractant in $L^1(\R^N)$. Hence, for all $t\geq 0$, $u(t)\in v + L^1(\R^N)$, and 
$$
\| u(t) - v\|_1\leq \| u_{ini} - v\|_1.
$$
Thus it is natural to compute the equation satisfied by $f(t)=u(t)-v\in L^1(\R^N)$: since $v$ is a stationary solution of \eqref{LCS}, there holds
$$
\p_t f + \dv_yB(y,f)- \Delta_y f=0,\quad t>0, y\in\R^N,
$$
where the flux $B$ is defined by
$$
B(y,f)=A(y,v(y)+ f ) - A(y,v(y)),\quad\forall (y,f)\in\T^N\times \R.
$$
The flux $B(y,f)$ vanishes at $f=0$, for all $f$. Moreover, if the flux $A$ satisfies the assumptions of Theorem \ref{thm:conv}, there exists $\alpha_1\in \mathcal C^1(\T^N)$ and $\tilde B_1\in \mathcal C(\T^N\times \R)$ such that
$$
   B(y,f)=\alpha_1(y) f + \tilde B_1(y,f),
$$
and the flux $\tilde B_1$ is such that
$$
\forall X>0, \ \exists C_X>0,\ \forall f\in[-X,X] ,\ \forall y\in\T^N,\ \left|\tilde B_1(y,f)\right|\leq C_X |f|^2.
$$
At some point in the proof, we will need a more refined approximation of $B$ in a neighbourhood of $f=0$; we thus also introduce $\alpha_2,\alpha_3\in L^\infty(\T^N)$, $\tilde B_3\in L^\infty(\T^N\times \R)$ such that
$$
B(y,f)=\alpha_1(y) f +\alpha_2(y)f^2 + \alpha_3(y) f^3 + \tilde B_3(y,f),
$$
and the flux $\tilde B_3$ is such that for all $X>0$, there exists a constant $C_X>0$, such that for all $f\in[-X,X]$, for all $y\in\T^N$, 
$$\begin{aligned}
 \left|\tilde B_3(y,f)\right|\leq C_X |f|^4,\\
\left|\dv_y\tilde B_3(y,f)\right|\leq C_X |f|^4,\\
 \left|\p_f\tilde B_3(y,f)\right|\leq C_X |f|^3.
  \end{aligned}
$$
The existence of $\alpha_i$ ($i=1,2,3$)  and the bounds on $\tilde B_1, \tilde B_3$ are ensured by the assumption that $A\in W^{5,\infty}(\T^N\times \R)$.
Notice in particular that 
$$
\alpha_1(y)=\p_f B(y,f)_{|f=0}= (\p_p A) (y,v(y)),\quad \forall y\in\T^N.
$$
\vskip1mm

As explained in \cite{BDK}, the interplay between the diffusion and the drift $\alpha_1$ induces a displacement of the center of mass. In the linear case, that is, when $\tilde B_1=0$, the evolution of the center of mass can be computed as follows: since the function $f$ satisfies
$$
\p_t f + \dv_y (\alpha_1 f )-\Delta_y f=0,
$$
there holds
$$
\frac{d}{dt}\int_{\R^N} y f(t,y)\:dy= N \int_{\R^N} \alpha_1(y) f(t,y)\:dy.
$$
Now, for $t\geq 0, y\in\T^N$, set 
$$
\tilde f(t,y)=\sum_{k\in\Z^N}f(t,y+k).
$$
Since the function $\alpha_1$ is periodic, $\tilde f$ satisfies
$$
\p_t \tilde f + \dv_y(\alpha_1 \tilde f) - \Delta_y \tilde f=0,\quad t>0, \ y\in \T^N,
$$
and we have, for all $t\geq 0$,
$$
\int_{\R^N}\alpha_1 f(t)=\int_{\T^N}\alpha_1 \tilde f(t).
$$

Using Lemma 1.1 of \cite{MMP} together with a Poincar\'e inequality on the torus $\T^N$, it can be easily proved that  as $t\to\infty$, $\tilde f(t)$ converges with exponential speed  in $L^1(\T^N)$ towards $\mean{\tilde f} f_0$, where $f_0$ is the unique solution of \eqref{eq:f0}.
 Additionally, notice that
$$
\mean{\tilde f}= \sum_{k\in\Z^N} \mean{ f(\cdot + k)}= \int_{\R^N} f=M.
$$

Consequently, setting
$$
c:= N\mean{\alpha_1 f_0}
$$
we infer that in the linear case,
$$
\frac{d}{dt}\int_{\R}(y-ct) f \to 0\quad \text{exponentially fast.}
$$

In fact, it turns out that the nonlinearity has no effect on this displacement, although this is not quite clear if we try to include the quadratic term $\tilde B_1$ in the above calculation. 
We will justify this result by formal calculations in the next paragraph.
Nonetheless, it can be proved in the case $N=1$ (see for instance \cite{longtime}) that when $\|f_0\|_1$ is not too large, 
$$
\| f(t)\|_{L^2(\R)}\leq C \frac{ \| f_0\|_1}{t^{1/4}}\quad\forall t>0,
$$  
and more generally, the $L^p$ norm of $f(t)$ vanishes for all $p\in(1,\infty]$.
This somehow explains why the quadratic term does not modify the motion of the center of mass for large times: the term $\tilde B_1(\cdot, f(t,\cdot))$ vanishes in $L^1(\R)$ as $t\to\infty$. Hence, hereinafter, we choose to make in the general case the same change of variables as the one dictated by the linear case. Precisely, let $U\in L^\infty_\text{loc}([0,\infty)\times \R^N)$ such that
\be\label{scaling}
f(t,y)= \frac{1}{(1+2t)^{N/2}} U\left(\log \sqrt{1+2t}, \frac{y-ct}{\sqrt{1+2t}}\right),\quad t\geq 0,\ y\in\R^N.
\ee
This change of variables is classical in the study of long-time parabolic dynamics, see for instance \cite{EZ}. In the present case, our change of variables is exactly the same as in \cite{BDK}; straightforward calculations lead to
\be\label{eq:rescaled}
\p_\tau U- \dv_x(x U) + R \dv_x((\alpha_1(z)-c) U) -\Delta_x U= - R^{N+1}\dv_x \tilde B_1\left( z,\frac{U}{R^N} \right),
\ee
with $\tau>0$, $x\in\R^N$, and where
$$
R= e^\tau\quad\text{and}\quad z=Rx+ c \frac{R^2-1}{2}.
$$
Studying the long time behaviour of $f$ amounts to studying the long time behaviour of $U$. Now, as $\tau\to \infty$, the quantity $R$ becomes very large, and thus the variable $z$ is highly oscillating. Hence, as emphasized in \cite{BDK}, the asymptotic study of equation \eqref{eq:rescaled} somehow falls into the scope of homogenization theory; the small parameter measuring the period of the oscillations is then $\e=R^{-1}=e^{-\tau}$. However, one substantial difference with classical homogenization problems is that the small parameter depends on time, which sometimes makes the proofs much more technical. We refer to \cite{BDK} for more details.

Let us also mention that the homogenization of equation \eqref{eq:rescaled} with a ``fixed'' small parameter, and when the quadratic flux $\tilde B_1$ vanishes, has been performed by Thierry Goudon and Fr\'ed\'eric Poupaud in \cite{GP}. As a consequence, the formal asymptotic expansions which will be performed in the next section are in fact very close to the ones of \cite{GP}.

\subsection{Formal derivation of the limit system}

As usual in homogenization problems (see \cite{BLP} for instance), the idea is now to assume that the solution $U$ of \eqref{eq:rescaled} admits an asymptotic development in powers of the small parameter measuring the period of the oscillations; in the present case, the small parameter is $e^{-\tau}$, so that we expect the approximation to be valid for large times only. Hence, assume that when $\tau\gg 1$,
\be
U(\tau,x)\approx U_0\left(\tau,x,z \right) + e^{-\tau} U_1\left(\tau,x,z \right) + e^{-2\tau} U_2\left(\tau,x,z \right)+\cdots \label{Ansatz}
\ee
where $z=e^\tau x + c (e^{2\tau} -1)/2$ stands for the fast  variable and where for all $(\tau,x)\in\R_+\times \R^N$, the function
$$
z\mapsto U_i(\tau,x,z)
$$
is $\T^N$-periodic. Plugging the Ansatz \eqref{Ansatz} into equation \eqref{eq:rescaled} and identifying the powers of $R= e^\tau$ leads to a cascade of equations on the terms $U_0, U_1$, etc. Notice that according to Lemma A.1 in the Appendix, $f\in L^\infty([0,\infty)\times \R^N)$, and thus $U/ R^N$ is bounded in $L^\infty$.

\vskip1mm
\noindent $\bullet$ \textbf{Terms of order $R^2$:} Identifying the highest order terms in equation \eqref{eq:rescaled} when $U$ is given by \eqref{Ansatz} leads to 
$$
c \cdot \nabla_z U_0 + \dv_z((\alpha_1-c) U_0)-\Delta_z U_0 = - \Delta_z U_0 + \dv(\alpha_1 U_0)=0,\quad z\in\T^N.
$$
We recall the following result, which is a straightforward consequence of the Krein-Rutman Theorem (see \cite{DL}):
\begin{lem}
Let $\alpha\in L^\infty(\T^N)$. Consider the vector space
$$
E[\alpha]:=\left\{w\in H^1(\T^N),\ -\Delta_z w + \dv_z (\alpha w)=0 \right\}.
$$
Then $\mathrm{dim} E[\alpha]=1$, and there exists a unique function $m\in E[\alpha]$ such that $\mean{ m}=1$.

Moreover, $m\in W^{1,p}(\T^N)$ for all $p<\infty$, and
$$ \inf_{z\in\T^N} m>0.$$
\end{lem}
In the present case, $E[\alpha_1]=\R f_0$, where $f_0$ is defined by \eqref{eq:f0}.
Hence there exists a function $F=F(\tau,x)$ such that
\be\label{def:U_0}
U_0(\tau,x,z)=f_0(z) F(\tau,x)\quad\forall (\tau,x,z)\in[0,\infty)\times \R^N\times \T^N.
\ee
\vskip1mm

\noindent $\bullet$ \textbf{Terms of order $R^1$:}
Concerning the terms of order $R^1= e^\tau$, the case when the space dimension is equal to one has to be treated separately. Indeed, 
\begin{eqnarray*}
	R^{N+1}\dv_x \tilde B_1\left(z, \frac{U}{R^N}  \right)&=& R^{1-N}\dv_x(\alpha_2 U^2 )+ R^{1-2N}\dv_x(\alpha_3 U^3) \\&&+ R^{N+1}\dv_x \tilde B_3\left(z, \frac{U}{R^N}\right),
\end{eqnarray*}
and using the bounds on $\tilde B_3$,
\begin{eqnarray*}
 R^{N+1}\dv_x \tilde B_3\left(z, \frac{U}{R^N}\right)&=& R^N (\dv_ z\tilde B_3 )\left(z, \frac{U}{R^N}\right) + R \nabla_x U \cdot (\p_U \tilde B_3 )\left(z, \frac{U}{R^N}\right) \\
&=& \mathcal O(R^{-3N}) + \mathcal O(R^{-3N+2})= \mathcal O(R^{-1}).
\end{eqnarray*}
We infer that if $U$ is given by \eqref{Ansatz}, 
\begin{eqnarray}
 \label{nonlinear_term}
R^{N+1}\dv_x \tilde B_1\left(z, \frac{U}{R^N}  \right)& =& R^{2-N}\dv_z(\alpha_2 U_0^2) \\&+& R^{1-N}\left[\dv_x(\alpha_2 U_0^2) + 2 \dv_z(\alpha_2 U_0 U_1)\right]\\\nonumber& + &R^{2-2N}\dv_z(\alpha_3 U_0^3)\nonumber\\&+&\nonumber  \mathcal O(R^{-1}).
\end{eqnarray}

Consequently, we obtain that when $N\geq 2$, the term $U_1$ solves the equation
\be\label{eq:U_1_Ngrand}
-\Delta_z U_1 + \dv_z(\alpha_1 U_1)= - \dv_x((\alpha_1-c ) U_0) + 2 \sum_{i=1}^N \frac{\partial^2 U_0}{\partial x_i \partial z_i}.
\ee
Since $U_0(t,x,z)=f_0(z) F(t,x) $, we have
$$
\mean{(\alpha_1 - c)U_0(t,x,\cdot)}= F(t,x)\left(\mean{\alpha_1 f_0} - c\right)=0
$$
by definition of $c$. Hence the right-hand side of \eqref{eq:U_1_Ngrand} has zero mean value, and the compatibility condition is satisfied. Thus for all $(t,x)\in[0,\infty)\times \R^N$, \eqref{eq:U_1_Ngrand} has a unique solution in $H^1(\T^N)$. Moreover, using the linearity of \eqref{eq:U_1_Ngrand} together with the expression \eqref{def:U_0}, we infer that $U_1$ can be written as
\be\label{def:U_1_Ngrand}
U_1(t,x,z)= f_1(z) \cdot \nabla_x F(t,x),
\ee
where $f_1\in H^1(\T^N)^N$ satisfies
\be\label{eq:f_1}
-\Delta_z f_{1,i} + \dv_z (\alpha_1f_{1,i})= - f_0 (\alpha_{1,i} - c_i) + 2 \p_{z_i} f_0,\quad\forall i\in\{1,\cdots, N\}.
\ee
Notice that according to the regularity assumptions on the flux $A$, the function $\alpha_1$ belongs to $ W^{1,\infty}(\T^N)$; thus $f_0\in W^{2,p}(\T^N)$ for all $p<\infty$, and therefore $f_1\in W^{2,p}(\T^N)$ for all $p<\infty$. In particular, $f_1, f_0\in W^{1,\infty}(\T^N).$

If $N=1$, on the other hand, the corrector $U_1$ solves the equation
\be
\label{eq:U_1_N=1}
-\p_{zz} U_1 + \p_z(\alpha_1 U_1)= - \p_x((\alpha_1-c ) U_0) + 2 \frac{\partial^2 U_0}{\partial x \partial z} - \p_z(\alpha_2 U_0^2).
\ee
Notice that the compatibility condition is satisfied, for the same reason as before.
Hence in this case,
\be
U_1(t,x,z)= f_1(z) \p_x F(t,x) + g_1(z) F(t,x)^2,\label{def:U_1_N=1}
\ee
where $g_1\in H^1(\T)$ solves
$$
-\Delta_z g_1 + \p_z(\alpha_1 g_1) = - \p_z(\alpha_2 (f_0)^2).
$$

The fact that the compatibility condition is satisfied in all cases justifies the use of the change of variables \eqref{scaling} in the nonlinear case. This means that, at least on a formal level, the displacement of the center of mass of the function $f$ is unaffected by the presence of the quadratic term $\tilde B_1$.

\vskip1mm
\noindent $\bullet$ \textbf{Terms of order $R^0$:}
As we identify the terms of order one in equation \eqref{eq:rescaled}, we obtain
\begin{eqnarray}
\label{eq:U2}&&- \Delta_z U_2 + \dv(\alpha_1 U_2)\\
\nonumber&=&-\p_\tau U_0 + \dv_x(xU_0)+ \Delta_x U_0 - \dv_x((\alpha_1-c) U_1) + 2 \sum_{i=1}^N \frac{\p^2 U_1}{\p x_i \p z_i} + \mathcal A_{NL},
\end{eqnarray}
where the term $\mathcal A_{NL}$ stems from the expansion of the nonlinear term $\tilde B_1$. According to \eqref{nonlinear_term}, we have
$$\begin{aligned}
   \mathcal A_{NL}=\p_x(\alpha_2 U_0^2) + 2 \p_z (\alpha_2U_0 U_1) + \p_z(\alpha_3 U_0^3)\quad \text{if }N=1,\\ 
\mathcal A_{NL}=\dv_z(\alpha_2 U_0^2)\quad \text{if }N=2,\\
\mathcal A_{NL}=0 \quad \text{if } N\geq 3.
  \end{aligned}
$$
The evolution equation  for the function $F$ follows from the compatibility condition; precisely, we obtain
$$
\p_\tau F - \dv_x(x F) - \Delta_x F + \dv_x\mean{(\alpha_1-c) U_1} - \mean{   \mathcal A_{NL}}=0.
$$
We now distinguish between the cases $N\geq 2$ and $N=1$.

$ \triangleright $ \underline{If $N\geq 2$}, $\mean{ \mathcal A_{NL}}=0$; using \eqref{def:U_1_Ngrand}, we infer that $F$ satisfies
\be\label{eq:h_Ngrand}
\p_\tau F - \dv_x(x F) - \sum_{1\leq i,j\leq N} \eta_{i,j}\frac{\p^2 F}{\p x_i \p x_j} =0,\quad \tau>0,\ x\in\R^N\text{ with }N\geq 2,
\ee
where the coefficients $(\eta_{i,j})_{1\leq i,j\leq N}$ are given by
$$
\eta_{i,j}= \delta_{i,j} - \mean{(\alpha_{1,i}-c_i) f_{1,j}}.
$$
The following Lemma entails that  equation \eqref{eq:h_Ngrand} is well-posed (see also Lemma \ref{prop:ex_uni_h}):
\begin{lem}
The matrix $\eta:=(\eta_{i,j})_{1\leq i,j\leq N}$ is coercive.
\label{lem:coercivity}
\end{lem}
Lemma \ref{lem:coercivity} is proved in \cite{GP} in dimension $N$, and its proof is recalled in \cite{BDK} when $N=1$. For the reader's convenience, we sketch the main steps of the proof here, and we refer to \cite{GP}, Proposition 4.6 for details.

\begin{proof}
Let $L$ be the differential operator
$$
L\phi=-\Delta_z \phi  +\dv_z (\alpha_1 \phi).
$$
The idea is to introduce, for all $j\in\{1,\cdots, N\}$, the function $\chi_j$ which solves the adjoint problem 
$$
L^* \chi_j=- \Delta_z \chi_j - \alpha_1\cdot \nabla_z \chi_j = \alpha_{1,j} - c_j,\quad \mean{\chi_j}=0.
$$
Since the right-hand side satisfies $\mean{(\alpha_{1,j} - c_j)\psi}=0$ for all $\psi\in\ker L=\R f_0$, 
the function $\chi_j$ is well-defined. For all $\xi\in \R^N$, we have
\begin{eqnarray*}
 \sum_{i,j}\mean{(\alpha_{1,i}-c_i) f_{1,j}}\xi_i\xi_j&=&\mean{L^*(\chi\cdot \xi) f_1\cdot \xi} \\&=& \mean{\chi\cdot \xi \left(- f_0 (\alpha_1-c)\cdot \xi  + 2 \nabla_z f_0 \cdot \xi \right)}\\
&=&-\mean{f_0\chi\cdot \xi L^*(\chi\cdot \xi)} + \mean{ 2 \chi\cdot \xi \nabla_z f_0 \cdot \xi}\\
&=&-\mean{L(f_0\chi\cdot \xi) \chi\cdot \xi} - 2 \mean{f_0  \xi \cdot\nabla_z (\chi\cdot \xi)}.
\end{eqnarray*}
Expanding $L(f_0\chi\cdot \xi) $ and using the identity $L f_0=0$ leads to
$$
\mean{L(f_0\chi\cdot \xi) \chi\cdot \xi}= \mean{f_0 \left| \nabla_z(\chi\cdot \xi) \right|^2}.
$$
Hence
\begin{eqnarray*}
	\sum_{1\leq i,j\leq N}\eta_{i,j}\xi_i \xi_j&=& |\xi|^2 + \mean{f_0 \left| \nabla_z(\chi\cdot \xi) \right|^2} + 2 \mean{f_0  \xi \cdot\nabla_z (\chi\cdot \xi)}\\&=&\mean{f_0 \left| \xi + \nabla_z(\chi\cdot\xi)\right|^2}.
\end{eqnarray*}
We deduce that
$$
 \sum_{1\leq i,j\leq N}\eta_{i,j}\xi_i \xi_j\geq 0\quad\forall \xi\in\R^N.
$$
Now, let $\xi\in\R^N$ such that $\sum \eta_{i,j}\xi_i \xi_j= 0$. Since $f_0(z)>0$ for all $z$, we infer that 
$$
\xi + \nabla_z(\chi\cdot\xi) =0\quad \forall z\in\T^N.
$$
Taking the average of the above inequality on $\T^N$ leads to $\xi=0$. Hence the matrix $(\eta_{i,j})$ is coercive.
\end{proof}

\vskip1mm

$\triangleright$ \underline{If $N=1$},  we have
$$
\mean{\mathcal A_{NL}}= \p_x \mean{\alpha_2 U_0^2} = \mean{\alpha_2 f_0^2}\p_x F^2.
$$
Moreover, in this case $U_1$ is given by \eqref{def:U_1_N=1}; hence
$$
 \mean{\dv_x((\alpha_1-c) U_1)}=\mean{(\alpha_1-c) f_1}\p_{xx} F + \mean{(\alpha_1-c) g_1} \p_xF^2.
$$
Consequently, the compatibility condition reads
\be\label{eq:h_N=1}
\p_\tau F - \p_x (xF) + a \p_x F^2- \eta \p_{xx} F=0,\quad \tau>0, \ x\in\R,
\ee
where the coefficients $a,\eta$ are given by
$$\begin{aligned}
a:=\mean{\alpha_2 f_0^2} + \mean{(\alpha_1-c) g_1},\\
\eta:=1 -\mean{(\alpha_1-c) f_1}.
\end{aligned}
$$
Lemma \ref{lem:coercivity} states that the diffusion coefficient $\eta$ is positive.

\vskip1mm

This completes the formal derivation of an approximate solution. In the following paragraphs, we recall or prove several results concerning the well-posedness and the long time behaviour of equations \eqref{eq:h_Ngrand} and \eqref{eq:h_N=1}. We will often refer to the equation on $F$ as the ``homogenized equation''; this term denotes equation \eqref{eq:h_Ngrand} when $N\geq 2$, and \eqref{eq:h_N=1} when $N=1$.

\subsection{Existence and uniqueness of stationary solutions}

This paragraph is concerned with the existence and uniqueness (in suitable functional spaces) of stationary solutions of the homogenized equations \eqref{eq:h_Ngrand} and \eqref{eq:h_N=1}. In the case when $N=1$, or when $(\eta_{i,j})_{1\leq i,j\leq N}=\lambda I$ for some $\lambda >0$, such results are stated in \cite{AEZ}. In the general case, we merely use a linear change of variables, and the problem is then reduced to the case of an isotropic diffusion.

\begin{lem}
Assume that $N\geq 2$. For $\gamma >0$, set $\psi_\gamma:x\in\R^N\mapsto \exp(\gamma |x|^2).$ Then there exists $\gamma>0$ such that for all $M\in\R$, there exists a unique function $F_M\in H^1(\psi_\gamma)$ satisfying
\be\label{eq:hM}
- \sum_{1\leq i,j\leq N} \eta_{i,j}\p_i \p_j F_M - \dv_x(xF_M)=0,\quad \int_{\R^N} F_M=M.
\ee
Furthermore, the following properties hold:
\begin{enumerate}[(i)]
 \item For all $M\in\R$, $F_M=Mh_1$;
\item $h_1\in W^{2,p}\cap \mathcal C^\infty(\R^N)$ for all $p\in[1,\infty)$, and $h_1\in H^2(\psi_\gamma);$
\item $h_1(x)>0$ for all $x\in\R^N$.
\end{enumerate}

\label{lem:statsol_Ngrand}
\end{lem}

\begin{proof}
The idea is to perform an affine change of variables in order to transform the diffusion term into a laplacian.  Precisely, set
$$
s_{i,j}= \frac{\eta_{i,j} + \eta_{j,i}}{2},\quad 1\leq i,j\leq N.
$$
Then the matrix $S=(s_{i,j})$ is symmetric and positive definite (see Lemma \ref{lem:coercivity}); hence there exists an orthogonal matrix $O\in\mathcal M_n(\R)$ and positive numbers $\lambda_j$ such that
$$
S=O^T \mathrm{Diag}(\lambda_1,\cdots, \lambda_N) O.
$$
Let us change the variables by setting
\be\label{CV}
x= P y, \quad\text{with }P:=O^T\mathrm{Diag}(\lambda_1^{1/2},\cdots, \lambda_N^{1/2})  ,
\ee
and for any function $F\in L^1(\R^N)$, define
$$
\tilde F(y)=F(Py).
$$
It can be readily checked that for all $x\in\R^N$,
$$
\sum_{i,j} \eta_{i,j}\frac{\p^2 F}{\p x_i \p x_j}(x)=\sum_{k,l} \tilde \eta_{k,l} \frac{\p^2 \tilde F }{\p y_k\p y_l} (P^{-1} x),
$$
where the coefficients $\tilde \eta_{k,l}$ are given by
$$
\tilde \eta_{k,l}=\sum_{i,j} (P^{-1})_{k,i} (P^{-1})_{l,j} s_{i,j} =(P^{-1} S (P^{-1})^T)_{k,l}.
$$
Using the definitions of the matrices $P$ and $S$, we infer that
$$
\tilde \eta=P^{-1} S (P^{-1})^T =I_N.
$$
Thus the diffusion term is transformed into a laplacian with this change of variables.

Let us now compute the drift term. We have 
$$
x_i= (Py)_i,
$$
and, denoting by $(e_1,\cdots, e_N)$ the canonical basis of $\R^N$,
\begin{eqnarray*}
\frac{\p F(x)}{\p x_i} &=&\frac{\p}{\p x_i}\tilde F(P^{-1} x)\\
&=&(P^{-1} e_i)\cdot \nabla_y \tilde F (P^{-1} x).
\end{eqnarray*}
Thus, always setting $x=Py$, 
\begin{eqnarray*}
x\cdot \nabla_x F(x)&=&\sum_{i=1}^N (Py)_i (P^{-1} e_i)\cdot \nabla_y \tilde F(y)\\
&=&\left[P^{-1} \left(\sum_{i=1}^N (P y)_i e_i\right)\right]\cdot \nabla_y \tilde F(y)\\
&=& (P^{-1} P y)\cdot \nabla_y \tilde F(y)=y\cdot \nabla_y \tilde F(y).
\end{eqnarray*}
Notice that this property is in fact independent of the definition of the matrix $P$.
Consequently, $F_M$ is a solution of \eqref{eq:hM} if and only if $\tilde F_M$ satisfies
$$
-\Delta_y \tilde F_M - \dv_y(y \tilde F_M)=0,\quad \int_{\R^N} \tilde F_M= (\det S)^{-1/2} M.
$$
The only solutions of the above equation in $H^1(\R^N)$ are Gaussian functions. Hence there exists a unique solution of \eqref{eq:hM} in $H^1(\R^N)$ for all $M$, and this solution is given by
$$
F_M(x)=C M \exp\left( \frac{1}{2}|P^{-1} x|^2\right),
$$
where the positive constant $C$ is a normalization factor. Moreover,
$$
|P^{-1}x|^2= \left|\mathrm{Diag}\left(\lambda_1^{-1/2}, \cdots, \lambda_N^{-1/2}\right)Ox \right|^2,
$$
and thus, since $|Ox|^2=|x|^2$,
$$
\left( \max_{1\leq i\leq N} \lambda_i\right)^{-1/2} |x|^2 \leq|P^{-1}x|^2\leq \left( \min_{1\leq i\leq N} \lambda_i\right)^{-1/2} |x|^2.
$$

All the properties of the lemma follow, with
$$
\gamma <2 \left( \max_{1\leq i\leq N} \lambda_i\right)^{-1/2}.
$$

\end{proof}

\vskip2mm

In the case when $N=1$, the existence of a stationary solution is treated in \cite{AEZ}. Hence we merely recall the main results of \cite{AEZ} in that regard.

\begin{lem}[Aguirre, Escobedo, Zuazua]
 Let $M\in\R$ be arbitrary, and let $a\in\R$, $\eta>0$. Let $\gamma:=(2\eta)^{-1}.$

Then there exists a unique function $F_M\in H^1(\psi_\gamma)$ which satisfies
$$
-\eta\p_{xx} F_M - \p_x(x F_M) + a \p_x F_M^2=0,\quad \int_{\R} F_M=M. 
$$
Moreover, $F_M$ enjoys the following properties:
\begin{enumerate}[(i)]
\item $F_M\in W^{2,p}\cap \mathcal C^\infty(\R)$ for all $p\in[1,\infty)$, and $F_M\in H^2(\psi_\gamma);$
\item If $M>0$, then $F_M(x)>0$ for all $x\in\R^N$.
\end{enumerate}
\label{lem:statsol_N=1}

\end{lem}

We deduce from the above Lemma that if $\gamma'<\gamma$, then there exists a constant $C_{\gamma'}$ such that 
$$
\left|F_M(x)  \right|,\ \left|\p_x F_M(x)  \right|\leq C_{\gamma'} \exp(-\gamma' x^2)\quad\forall x\in\R.
$$
Indeed, since $F_M\in H^2(\psi_\gamma)$, it can be easily proved that $F_M\psi_{\gamma'}\in H^2(\R)$ for all $\gamma'<\gamma$. Sobolev embeddings then imply that $F_M\psi_{\gamma'}\in W^{1,\infty}(\R).$

The existence of stationary solutions of \eqref{eq:h_Ngrand} and \eqref{eq:h_N=1} is now ensured. We now tackle the study of the properties of equations \eqref{eq:h_Ngrand} and \eqref{eq:h_N=1}, focusing in particular on the long-time behaviour and on regularity issues.

We begin with a definition of the weight function $K\in \mathcal C^\infty(\R^N)$, which plays a central role in the theory of existence. We use  the change of variables \eqref{CV}, which was introduced in the proof of Lemma \ref{lem:statsol_Ngrand}. This allows us to transform the matrix $(\eta_{i,j})$ into the identity matrix. 
For $t>0$, $y\in\R^N$, set 
$\tilde F(t,y)=F(t,P y)$. If $F$ is a solution of \eqref{eq:h_Ngrand}, then $\tilde F$ solves
$$
\p_t \tilde F -\dv_y(y\tilde F) -\Delta_y \tilde F=0.
$$
Consequently, the results of \cite{EZ} can be directly applied to $\tilde F$, for which existence is proved in the functional space $L^2(K_0)$, where $K_0(y)=\exp(y^2/2)$.
 Performing the inverse change of variables, it is clear that the relevant weight function is given by
$$
K(x):=K_0(P^{-1}x)=\exp\left(\frac{|P^{-1} x|^2}{2}\right).
$$
Notice that by definition of the matrix $P$, there exist positive constants $\gamma, \gamma'$ such that
$$
\exp(\gamma' x^2)\leq K(x)\leq \exp(\gamma x^2)\quad\forall x\in\R^N.
$$
When $N=1$, the weight function $K$ is given by
$$
K(x):=\exp\left(\frac{|x|^2}{2\eta}\right).
$$

We immediately deduce from \cite{EZ} the following Proposition:
\begin{prop}
Let $F_{ini}\in L^\infty(\R^N)\cap L^2(K)$. Then the homogenized problem has a unique solution
$$
F\in \mathcal C([0,\infty), L^2(K))\cap \mathcal C((0,\infty), H^2(K))\cap \mathcal C^1((0,\infty), L^2(K))
$$
such that $F_{|t=0}=F_{ini}$.

Moreover, 
$$
\lim_{t\to \infty} \|F(t) - F_M \|_{L^1(\R^N)}=0,
$$
where $F_M$ is the unique stationary solution of the homogenized problem with mass $M=\int_{\R^N} F_{ini}$.

\label{prop:ex_uni_h}
\end{prop}

Consequently, the homogenized equations \eqref{eq:h_N=1} and \eqref{eq:h_Ngrand} are well posed. We conclude this section by stating a result on the construction of an approximate solution:
\begin{defi}
Let $F\in \mathcal C([0,\infty, L^2(K))\cap \mathcal C((0,\infty), H^2(K)).$ We define the approximate solution of \eqref{eq:rescaled} associated with $F$ by
$$
\Uapp[F](\tau,x; R)=U_0\left( \tau, x, z \right) + R^{-1} U_1\left( \tau,x, z \right)  + R^{-2} U_2\left(\tau, x, z \right),
$$
with $\tau\geq 0$, $x\in\R^N$, $R>0$ and $z:=Rx + c\frac{R^2-1}{2}$, and where 
\begin{itemize}
 \item $U_0$ is defined by \eqref{def:U_0};
\item $U_1$ is defined by \eqref{def:U_1_Ngrand} if $N\geq 2$ and by \eqref{def:U_1_N=1} if $N=1$;
\item $U_2$ is defined by 
\begin{eqnarray*}
-\Delta_z U_2 + \dv(\alpha_1 U_2) 
&=&(f_0(z)-1)\left[ -\p_\tau F + \dv_x(x F) +\Delta_x F \right]\\
&&+\mean{\dv_x((\alpha_1-c )U_1)} - \dv_x((\alpha_1-c )U_1)\\
&&+ 2 \sum_{i=1}^N\frac{\p^2 U_1}{\p x_i\p z_i} + \cA_{NL}-\mean{\cA_{NL}}.
\end{eqnarray*}

\end{itemize}
\label{def:Uapp}
\end{defi}
Notice that we do not require, in the above definition, that $F$ is a solution of \eqref{eq:h_Ngrand} or \eqref{eq:h_N=1}; hence the right-hand side in the equation on $U_2$ is slightly modified, so that the compatibility condition is satisfied and $U_2$ is well-defined. Of course, if $F$ is a solution of \eqref{eq:h_N=1} or \eqref{eq:h_Ngrand}, the equation on $U_2$ becomes \eqref{eq:U2}.

We then have the following result:
\begin{lem}
\begin{enumerate}
 \item Let $M\in\R$ be arbitrary. Define the function $U\in L^\infty([0,\infty)\cap \R^N))$ by $$U(\tau,x):=\Uapp[F_M]\left(\tau,x; e^{\tau}  \right).$$ 

Then $U$ is a solution of
\begin{multline*}
\p_\tau U -\dv_x(xU)-\Delta_x U +R \dv_x((\alpha_1(z)-c )U)=\\
=-R^{N+1}\dv_x \tilde B_1\left(z,\frac{U}{R^N}\right) + U^\text{rem},	
\end{multline*}
where the remainder term $U^\text{rem}$ is such that there exist $C>0$, $\gamma>0$ such that
$$
\left\|U^\text{rem}(\tau) \right\|_{L^1(e^{\gamma |x|^2})} + \left\|U^\text{rem}(\tau) \right\|_{L^\infty(\R^N)}\leq C e^{-\tau} \quad\forall \tau\geq 0.
$$
\item Let $F_{ini}\in L^\infty(\R^N)\cap L^2(K)$, and let $F\in\mathcal C([0,\infty), L^2(K))$ be the unique solution of the homogenized equation such that $F_{|t=0}=F_{ini}.$ Let $\rho\in\mathcal C^\infty_0(\R^N)$ be a mollifying kernel ($\rho\geq 0$, $\int \rho=1$), and let $F_\delta:=F\ast_x\rho_\delta$, where $\rho_\delta=\delta^{-N} \rho(\cdot/\delta)$, for $\delta>0$.

Let $(\tau_n)_{n\geq 0}$ be a sequence of positive numbers such that $\ds\lim_{n\to\infty} \tau_n=+\infty$. For $n\in\N, \delta>0$, define the function $u_n^\delta$ by
$$
u_n^\delta(\tau,x)=\Uapp[F_\delta]\left(\tau,x; e^{\tau_n+\tau} \right),\quad x\in\R^N, \ \tau\geq0.
$$
Then $u_n^\delta$ satisfies, with $R_n=e^{\tau_n+\tau}$ and $z_n=R_n x + c \frac{R_n^2-1}2$,
\begin{multline*}
	\p_\tau u_n^\delta -\dv_x(xu_n^\delta) +R_n \dv_x((\alpha_1(z_n)-c )u_n^\delta)-\Delta_x u_n^\delta=\\=-R_n^{N+1}\dv_x \tilde B_1\left(z,\frac{u_n^\delta}{R_n^N}\right) + r_n^\delta,
\end{multline*}
where the remainder term $r_n^\delta$ satisfies, for all $T>0$, 
$$
\| r_n^\delta\|_{L^\infty([0,T], L^1(\R^N))}\leq \omega_T(\delta) + C_{\delta,T} e^{-\tau_n},
$$
where $\omega_T:\R_+\to\R_+$ is a function depending only on $T$ such that $\lim_{0^+} \omega_T=0.$
\end{enumerate}

\label{lem:approx_sol}

\end{lem}
The proof of the above Lemma follows the calculations of the first paragraph; the proof is lengthy but straightforward, and is therefore left to the reader. The fact that $U^\text{rem}$ has exponential decay is a consequence of Lemmas \ref{lem:statsol_Ngrand}, \ref{lem:statsol_N=1}.

\section{Weighted $L^2 $ bounds for the rescaled equation}

\label{sec:bounds}

As explained in the previous section, we choose to work with the rescaled equation \eqref{eq:rescaled} rather than with the original one \eqref{LCS}. In fact, it can be easily checked that Theorem \ref{thm:conv} is equivalent to the following Proposition:

\begin{prop}
Let $U_{ini}\in L^1(\R^N)$, and let $M:=\int_{\R^N} U_{ini}$. 

Let $U\in \mathcal C([0,\infty), L^1(\R^N))$ be the unique solution of \eqref{eq:rescaled} with initial data $U_{|\tau=0}=U_{ini}.$ Then
$$
\lim_{\tau\to\infty}\int_{\R^N} \left| U(\tau,x) - f_0\left( e^{\tau} x + c \frac{e^{2\tau} -1}{2}  \right) F_M(x)\right|\:dx=0,
$$
 where the speed $c$ is defined by \eqref{def:c}, and $F_M\in L^1(\R^N)$ is the unique stationary solution of the homogenized equation \eqref{eq:h_Ngrand}, \eqref{eq:h_N=1} with total mass $M$.

\label{prop:conv_rescaled}
\end{prop}

In turn, since the function $f_0\in L^\infty(\T^N)$ is such that $\inf_{\T^N} f_0 >0$, the above statement is equivalent to
$$
\lim_{\tau\to\infty} \| V(\tau)-F_M\|_{L^1(\R^N)}=0,
$$
where the function $V=V(\tau,x)$ is defined by
\be\label{def:V}
V(\tau, x):=\frac{U(\tau,x) }{f_0\left( e^{\tau} x + c \frac{e^{2\tau} -1}{2}   \right)},\quad \tau\geq 0, \ x\in\R^N.
\ee

The proof of Proposition \ref{prop:conv_rescaled} consists of essentially two steps: first, we prove compactness properties in $L^1(\R^N)$ for  the family $(V(\tau))_{\tau \geq 0}$. To that end, we derive uniform bounds with respect to $\tau$ in weighted $L^2$ spaces; this step will be achieved in the current section. Then, we prove in the next section, using techniques inherited from dynamical systems theory, that the limit of any converging sequence $V(\tau_n)$ is equal to $F_M$. As emphasized in the introduction, the proof of convergence relies on rather abstract arguments, and thus does not yield any rate of convergence in general. However, when the flux $A$ is linear, the weighted $L^2$ bounds allow us to prove that the family $U(\tau)$ has uniformly bounded moments of order four, and thus \eqref{hyp:moments} holds. As proved in \cite{BDK}, the convergence stated in Theorem~\ref{thm:conv} then takes place with algebraic rate.

\vskip2mm

The main result of this section is the following:
\begin{prop}
 Let $U_{ini}\in L^1\cap L^\infty(\R^N)$, and let  $U\in \mathcal C([0,\infty), L^1(\R^N))$ be the unique solution of \eqref{eq:rescaled} with initial data $U_{|\tau=0}=U_{ini}.$ 
Let $m>2(N+1)$ be arbitrary, and assume that
$$
\int_{\R^N} |U_{ini}(x)|^2 (1+ |x|^2)^{m/2}\:dx<+ \infty.
$$
Then there exists a constant $C_m>0$ (depending only on $m, N,$ and on the flux $A$) such that if $\|U_{ini}\|_{L^1}\leq C_m,$ then
\be
\begin{aligned}
\sup_{\tau \geq 0} \int_{\R^N} |V(\tau,x)|^2(1+ |x|^2)^{m/2}\: dx<+ \infty,\\
\sup_{\tau\geq 0} \int_\tau^{\tau+1} \int_{\R^N} |\nabla_x V(s,x)|^2\:dx\:ds<+\infty.
\end{aligned}
\label{in:L2control}\ee
As a consequence, there exists a sequence $(\tau_n)$ of positive numbers  such that $\tau_n\in[n,n+1]$ for all $n$, and such that the sequence $(V(\tau_n,x))_{n\geq 0}$ is compact in $L^1(\R^N)$.

Moreover, if the flux $A$ is linear, then $C_m=+\infty$ for all $m>2(N+1)$.
\label{prop:compactness}
\end{prop}

Before proving the bounds \eqref{in:L2control}, we explain how they entail the existence of a converging sequence. Thus we admit that \eqref{in:L2control} holds for the time being. First, for any $X\geq 1$, $\tau\geq 0$, we have
\begin{eqnarray*}
&&\int_{|x|\geq X}|V(\tau,x)|\:dx\\&\leq & \left( \int_{|x|\geq X} |V(\tau,x)|^2(1+ |x|^2)^{m/2}\: dx \right)^{1/2} \left(\int_{|x|\geq X} (1+|x|^2)^{-m/2} \:dx \right)^{1/2}\\
&\leq & C X^{(N-m)/2} \left(\sup_{\tau \geq 0} \int_{\R^N} |V(\tau,x)|^2(1+ |x|^2)^{m/2}\: dx \right)^{1/2}.
\end{eqnarray*}
Since $m>N$, we infer that the family $\{ V(\tau,x)\}_{\tau\geq 0}$ is equi-integrable.

Moreover, let $\mathcal K\subset \R^N$ be an arbitrary compact set, and let $h\in\R^N$ be arbitrary, with $|h|\leq 1$. Let 
$$
\tilde{\mathcal K}:=\{x \in \R^N,\ d(x,\mathcal K)\leq 1\}.
$$
The set $\tilde{\mathcal K}$ is clearly compact. Then 
\begin{eqnarray*}
&&\int_{\mathcal K} |V(\tau,x+h) -V(\tau,x)|\:dx\\ &\leq & |h| \int_{\mathcal K}\int_0^1 |\nabla V|(\tau,x+\lambda h)\:d\lambda\:dx\\
&\leq & |h| \int_{\tilde{\mathcal K}} |\nabla V(\tau,z)|\:dz\leq |h| |\tilde{\mathcal K}|^{1/2} \left(\int_{\R^N} |\nabla_x V(\tau,x)|^2\:dx\right)^{1/2}.
\end{eqnarray*}

Now, for all $n\in\N$, there exists $\tau_n\in[n,n+1]$ such that 
$$
\int_{\R^N} |\nabla_x V(\tau_n,x)|^2\:dx \leq  \int_n^{n+1} \int_{\R^N} |\nabla_x V(s,x)|^2\:dx\:ds.
$$ 
Consequently, there exists a constant $C$, depending only on $\mathcal K$ and on the  bounds on $V$ in $L^2_\text{loc}([0,\infty), H^1)$, such that 
$$
\forall n\in\N,\quad \int_{\mathcal K} |V(\tau_n,x+h) -V(\tau_n,x)|\:dx \leq C |h|.
$$
Hence the sequence $\{ V(\tau_n,x)\}_{n\geq 0}$ is equi-continuous in $L^1(\R^N)$.

Notice also that
$$
\sup_{n\geq 0}\| V(\tau_n)\|_{L^1}\leq \frac{1}{\inf_{\T^N} f_0}\sup_{n\geq 0}\| U(\tau_n)\|_{L^1}\leq \frac{\|U_{ini}\|_1}{\inf_{\T^N} f_0}.
$$
Thus the sequence  $\{ V(\tau_n,x)\}_{n\geq 0}$ is bounded in $L^1(\R^N)$.

According to classical results of functional analysis (see for instance \cite{Brezis}), we infer that the sequence $(V(\tau_n))_{n\geq0}$ is compact in $L^1$.

\vskip1mm

The rest of the section is devoted to the proof of the bounds \eqref{in:L2control}. We first prove that $V\in L^\infty_\text{loc}([0,\infty), L^2((1+ |x|^2)^{m/2})).$ Then, using the construction of approximate solutions of \eqref{eq:rescaled} performed in the previous section, we derive an energy inequality on the function $V$. Carefully controlling the non-linear terms appearing in this energy inequality, we are led to \eqref{in:L2control}. 

Before addressing the proof, we recall a result which will play a key role in several arguments: since $U_{ini}\in L^\infty\cap L^1(\R^N)$, there exists a positive constant $C$, depending only on the flux $A$ and on $\| U_{ini}\|_1, \|U_{ini}\|_\infty$, such that
\be\label{in:Linfty_bound}
\|U(\tau)\|_{L^\infty(\R^N)}\leq C e^{N\tau}.
\ee
Indeed, performing backwards the parabolic scaling \eqref{scaling}, it turns out that this inequality is equivalent to the boundedness of $u$ in $L^\infty([0,\infty)\times \R^N)$, where $u$ is the solution of \eqref{LCS} with initial data $v + U_{ini}$. And the $L^\infty$ bound on $u$ follows from Lemma A.1 in the Appendix.

\vskip2mm

\noindent \textbf{First step: the family $V(\tau)$ is locally bounded in $L^2((1+ |x|^2)^{m/2})$.}

This amounts in fact to proving that $U\in L^\infty_\text{loc}([0,\infty), L^2((1+ |x|^2)^{m/2}))$. Hence, multiply \eqref{eq:rescaled} by $U(\tau,x) (1+|x|^2)^{m/2}$ and integrate with respect to the variable $x$. Always with the notation $R= e^\tau$, $z=R x + c \frac{R^2-1}{2}$, this leads to
\begin{eqnarray}\label{in:loc_bounds}
&&\frac{1}{2}\frac{d}{d\tau}\int_{\R^N} |U(\tau,x)|^2  (1+ |x|^2)^{m/2}dx\\
\nonumber&=&-\int_{\R^N} |\nabla_x U|^2 (1+ |x|^2)^{m/2}\:dx - {m}\int_{\R^N} (x\cdot \nabla_x U) U (1+ |x|^2)^{-1+m/2}\:dx\\
\nonumber&&-\frac{1}{2}\int_{\R^N}  (x\cdot \nabla_x |U|^2)(1+ |x|^2)^{m/2}\:dx - m\int_{\R^N}|U|^2 |x|^2(1+ |x|^2)^{-1+m/2}\:dx\\
\nonumber&&- R^{N+1}\int_{\R^N}\left[B\left(z , \frac{U(\tau,x)}{R^N } \right)  - c \frac{U(\tau,x)}{R^N }\right]\cdot \nabla_x  U (1+ |x|^2)^{m/2}dy\\
\nonumber&&- m  R^{N+1}\int_{\R^N}\left[B\left(z , \frac{U(\tau,x)}{R^N } \right)  - c \frac{U(\tau,x)}{R^N }\right]\cdot x \;U (1+ |x|^2)^{-1+m/2}\:dx.
\end{eqnarray}

Since $U(\tau)/R^N$ is bounded (see \eqref{in:Linfty_bound}), there exists a constant $C$ such that
$$
\left| B\left(Rx + c \frac{R^2-1}{2} , \frac{U(\tau,x)}{R^N } \right)  - c \frac{U(\tau,x)}{R^N }\right|\leq C \frac{|U(\tau,x)|}{R^N}.
$$
Moreover, 
$$
\left|x(1+|x|^2)^{-1+m/2}  \right|, \ \left||x|^2(1+|x|^2)^{-1+m/2}  \right|\leq (1+|x|^2)^{m/2}\quad\forall x\in\R^N.
$$
Hence, using the Cauchy-Schwarz inequality, we infer that the last two terms in \eqref{in:loc_bounds} are bounded by
$$
\frac{1}{4}\int_{\R^N} |\nabla_x U|^2 (1+ |x|^2)^{m/2}\:dx + C R^2 \int_{\R^N} | U|^2 (1+ |x|^2)^{m/2}\:dx.
$$
On the other hand,
\begin{eqnarray*}
&&\left| \int_{\R^N}  \left(x\cdot \nabla_x |U|^2\right)(1+ |x|^2)^{m/2}\:dx\right|\\ &=&\left|\int_{\R^N} |U|^2  \left( N(1+|x|^2)^{m/2} + m x^2(1+ |x|^2)^{-1+m/2}\right)\:dx\right|\\
&\leq &C \int_{\R^N} | U|^2 (1+ |x|^2)^{m/2}\:dx.
\end{eqnarray*}
Gathering all the terms, we obtain
\begin{eqnarray*}
&& \frac{d}{d\tau}\int_{\R^N} |U(\tau,x)|^2  (1+ |x|^2)^{m/2}dx\\
&\leq &- \int_{\R^N} |\nabla_x U(\tau,x)|^2 (1+ |x|^2)^{m/2}\:dx + C e^{2\tau}\int_{\R^N} |U(\tau,x)|^2  (1+ |x|^2)^{m/2}dx.
\end{eqnarray*}
Using Gronwall's Lemma, we deduce that 
$$
 U\in L^\infty_\text{loc}([0,\infty), L^2((1+|x|^2)^{m/2}),\quad \nabla_x U\in L^2_\text{loc}([0,\infty), L^2((1+|x|^2)^{m/2})).
$$

\vskip2mm
\noindent \textbf{Second step: The energy inequality.}
\vskip1MM

The idea here is the following: assume momentarily that the flux $B$ is linear, that is, $\tilde B_1=0.$ Let $\psi\in L^\infty([0,\infty)\times \R^N)$ be a solution of \eqref{eq:rescaled} such that $\psi(\tau,x)>0$ for all $\tau,x$. Then, according to \cite{MMP}, for any convex function $H\in\mathcal C^2(\R)$, we have
$$
\frac{d}{dt}\int_{\R^N} \psi(\tau,x) H\left(\frac{U(\tau,x)}{\psi(\tau,x)}\right)\:dx = - \int_{\R^N} H''\left(\frac{U(\tau,x)}{\psi(\tau,x)}\right)\left| \nabla_x \frac{U(\tau,x)}{\psi(\tau,x)}\right|^2\:dx.
$$
Taking $H:x\in\R\mapsto x^2$, we infer that
$$\sup_{\tau\geq 0}  \int_{\R^N} |U(\tau,x)|^2 \frac{dx}{\psi(x)}<+\infty.$$
Hence, if $\psi(x)$ behaves like $(1+|x|^2)^{-m/2}$ for $|x|$ large, the $L^2$ bound in  \eqref{in:L2control} is proved.

Thus, the goal of this step is to build a positive function $\tilde  U$, which behaves like  $(1+|x|^2)^{-m/2}$ for $|x|$ large, and which is an approximate solution of the linear part of \eqref{eq:rescaled}, with remainder terms of order one. Using calculations similar to the ones led in \cite{MMP}, we then derive an inequality on the energy
$$
\int_{\R^N} \left|\frac{U(\tau,x)}{\tilde U(\tau,x)}  \right|^2 \tilde U(\tau,x).
$$ 
From now on, we no longer assume that $\tilde B_1=0$.

The definition of $\tU$ is inspired from the construction of an approximate solution in the previous paragraph. Precisely, we set
$$
\tU (\tau,x)= f_0(z) h_m(x) + e^{-\tau}f_1(z)\cdot \nabla_x h_m(x),\quad \tau\geq 0, \ x\in\R^N, \ z=e^{\tau} x + c \frac{e^{2\tau}-1}{2},
$$
where the function $f_1\in W^{1,\infty}(\T^N)^N$ is defined by \eqref{eq:f_1}, and where
$$
h_m(x):=(1+ |x|^2)^{-m/2}.
$$
Remember that $\inf_{\T^N} f_0>0$; since
$$
\nabla_x h_m(x)= - m \frac{x}{1+ |x|^2} h_m,\quad y\in\R^N,
$$
we deduce that there exists $\tau_0>0$ (depending on $m$), such that 
\be\label{in:tU}
0<\frac{1}{2} f_0(z) h_m(x)\leq  \tU(\tau,x)\leq 2 f_0(z) h_m(x)\quad \forall y\in\R^N, \ \tau\geq \tau_0.
\ee
We now compute, for $\tau\geq \tau_0$, the rate of growth (or decay) of the energy $\int |U|^2 \tU^{-1}.$ Using equation \eqref{eq:rescaled} and performing several integrations by parts, we obtain
\begin{eqnarray*}
&&\frac{d}{d\tau}\int_{\R^N} \left|\frac{U(\tau,x)}{\tU(\tau,x)}  \right|^2 \tU(\tau,x)\:dx\\&=&-2\int_{\R^N} \left|\nabla_x\frac{U(\tau,x)}{\tU(\tau,x)}  \right|^2 \tU(\tau,x)\:dx\\
&&+\int_{\R^N} \left|\frac{U(\tau,x)}{\tU(\tau,x)}  \right|^2 \left[-\p_\tau \tU + \Delta_x \tU + \dv_x(x\tU) - e^{\tau}\dv_x\left((\alpha_1(z)-c)\tU\right)\right]\:dx\\
&&+2e^{(N+1)\tau}\int_{\R^N} \tilde B_1\left(z,\frac{U(\tau,x)}{e^{N\tau}}\right)\cdot \nabla_x \frac{U(\tau,x)}{\tU(\tau,x)}\:dx.
\end{eqnarray*}
By definition of $\tU$, we have
\begin{eqnarray*}
 &&-\p_\tau \tU + \Delta_x \tU + \dv_x(x\tU) - e^{\tau}\dv_x\left((\alpha_1(z)-c)\tU\right)\\
&=&\dv_x(xh_m) f_0(z)+f_0(z)\Delta_x h_m(x) + 2 \sum_{1\leq i,j\leq N} \frac{\p f_{1,i}}{\p z_j}(z)\frac{\p^2 h_m(x) }{\p x_i\p x_j}\\
&&-\sum_{1\leq i,j\leq N}\left[(\alpha_{1,i}-c_i) f_{1,j} \right](z)\frac{\p^2 h_m(x) }{\p x_i\p x_j}\\
&&+e^{-\tau}\sum_{1\leq i,j\leq N}f_{1,i}(z)\left[\frac{\p }{\p x_j}\left(x_j \p_{x_i} h_m(x)\right)+ \frac{\p^3 h_m(x)}{\p x_i \p x_j^2}\right],
\end{eqnarray*}
where
$$
z= e^\tau x + c \frac{e^{2\tau} -1}{2}.
$$
Notice that
$$\dv_x(xh_m(x))=(N-m) h_m(x) + \frac{m}{(1+|x|^2)^{1+\frac{m}{2}}},
$$

and there exists a constant $C$ (depending on $m$ and $N$) such that for all $i,j\in\{1,\cdots, N\}$,
$$\begin{aligned}
\left| \frac{\p^2h_m(x)}{\p x_i \p x_j}\right|\leq C\frac{1}{(1+ |x|^2)^{1+\frac{m}{2}}},\\
\left|  \nabla_x h_m(x)\right| + \left|  |x|\frac{\p^2h_m(x)}{\p x_i \p x_j}\right|+\left| \frac{\p^3 h_m(x)}{\p x_i \p x_j^2}\right|  \leq C h_m(x).
\end{aligned}
$$
Hence we infer (remember that $N-m<0$ and that inequality \eqref{in:tU} holds)
\begin{eqnarray*}
  &&-\p_\tau \tU + \Delta_x \tU + \dv_x(x\tU) - e^{\tau}\dv_x\left((\alpha_1(z)-c)\tU\right)\\
&\leq & (N-m) f_0(z) h_m(x)+ C e^{-\tau}h_m(x) + C\frac{1}{(1+ |x|^2)^{1+\frac{m}{2}}}\\
&\leq & \frac{N-m}{4}\tU(\tau,x) + C\frac{1}{(1+ |x|^2)^{1+\frac{m}{2}}}
\end{eqnarray*}
for all $\tau\geq \tau_0$, provided $\tau_0$ is chosen sufficiently large.

On the other hand, since the flux $\tilde B_1$ is quadratic near the origin and $U/e^{N\tau}$ is bounded, we have
\begin{eqnarray*}
&&	\left| e^{(N+1)\tau}\int_{\R^N} \tilde B_1\left(z,\frac{U(\tau,x)}{e^{N\tau}}\right)\cdot \nabla_x \frac{U(\tau,x)}{\tU(\tau,x)}\:dx\right|\\
&\leq &C e^{(1-N)\tau}\int_{\R^N} |U(\tau,x)|^2 \left|  \nabla_x \frac{U(\tau,x)}{\tU(\tau,x)} \right|\:dx.
\end{eqnarray*}
Gathering all the terms, we obtain
\begin{eqnarray}
&&\frac{d}{d\tau } \int_{\R^N} \left|\frac{U}{\tU}  \right|^2 \tU+ \frac{m-N}{4}\int_{\R^N} \left|\frac{U}{\tU}  \right|^2 \tU + 2 \int_{\R^N} \left|\nabla\frac{U}{\tU}  \right|^2 \tU\label{energy}\\
&\leq & C \int_{\R^N} \left(\frac{U(\cdot,x)}{\tU(\cdot,x)}\right)^2 \frac{dx}{(1+|x|^2)^{1+\frac{m}{2}}}\label{energy1} \\
&&+  C e^{(1-N)\tau}\int_{\R^N} |U(\cdot,x)|^2 \left|  \nabla_x \frac{U(\cdot,x)}{\tU(\cdot,x)} \right|\:dx.\label{energy2}
\end{eqnarray}
Notice that when the flux $A$ is linear, the term \eqref{energy2} is zero.

\vskip2mm
\noindent \textbf{Third step: control of the term \eqref{energy1}.}
\vskip1MM

Set $\phi:=U/\tU;$ then according to the first step, $\phi\in L^\infty_\text{loc}([\tau_0,\infty), L^2(h_m))\cap L^2_\text{loc}([\tau_0,\infty), H^1(h_m)).$ Moreover, 
$$
\nabla(\phi^2 h_m)= 2 \phi h_m \nabla \phi + \phi^2 \nabla h_m;
$$
since $|\nabla h_m|\leq m h_m,$ we deduce that $\phi^2 h_m\in L^1_\text{loc}([\tau_0,\infty), W^{1,1}(\R^N))$, and thus, using Sobolev embeddings, $\phi^2 h_m \in L^1_\text{loc}([\tau_0,\infty), L^{p^*}(\R^N)), $ where $p^*:=N/(N-1)$ if $N\geq 2$, and $p^*=+\infty$ if $N=1$. Additionally, the following inequality holds: there exists a constant $C$, depending only on $N$, such that for all  $\tau\geq \tau_0$
\begin{eqnarray*}
\| \phi^2 h_m(\tau)\|_{L^{p^*}(\R^N)}&\leq & C \|\nabla(\phi^2 h_m(\tau)) \|_{L^1(\R^N)}\\
&\leq & C \| \phi(\tau)\|_{L^2(h_m)}\| \nabla \phi(\tau)\|_{L^2(h_m)} + C \| \phi(\tau)\|_{L^2(h_m)}^2.
\end{eqnarray*}

We use the above inequality in order to control the term \eqref{energy1}.
First, let us write
\begin{eqnarray*}
 && \int_{\R^N} |\phi(\tau,x)|^2 (1+|x|^2)^{-\left(1+\frac{m}{2}\right)}\:dx\\
&=&\int_{\R^N} \left(\phi^2(\tau) h_m\right)^a \left( |\phi|(\tau) h_m \right)^b,
\end{eqnarray*}
where the exponents $a,b$ satisfy
$$
\left\{ \begin{array}{c}
         2a+b=2\\
	\ds a\frac{m}{2} + b \frac{m}{2} =1+\frac{m}{2},
        \end{array}
\right.$$
which leads to $a=1-\frac{2}{m}$, $b=\frac{4}{m}.$ Notice that $a, b\in (0,1)$, provided $m$ is large enough ($m>4$, which is always satisfied if $m>2(N+1)$).

Then, using H\"older's inequality, we infer
\begin{eqnarray*}
 && \int_{\R^N} |\phi(\tau,x)|^2 (1+|x|^2)^{-\left(1+\frac{m}{2}\right)}\:dx\\
&\leq &\left\|\phi^2(\tau) h_m \right\|_{L^p(\R^N)}^a \|  \phi(\tau) h_m \|_{L^1(\R^N)}^b,
\end{eqnarray*}
where the parameter $p$ is given by 
$$
p= a \left( 1 - b\right)^{-1}= \frac{1- \frac{2}{m}}{1-\frac{4}{m}}.
$$
Notice that $p$ is always larger than one. In order to be able to interpolate $L^p$ between $L^1$ and $L^{p^*}$, $p$ must also be smaller than $p^*$; if $N=1$, $p^*=\infty$, and thus we always have $p< p^*$. If $N\geq 2$, this condition amounts to $m> 2(N+1)$; we assume that $m$ always satisfies this assumption in the sequel.

Now, let $\theta\in(0,1)$ such that
$$
\frac{1}{p}=\frac{\theta}{1}+ \frac{1-\theta}{p^*};
$$
using once again H\"older's inequality, we obtain
\begin{eqnarray*}
  && \int_{\R^N} |\phi(\cdot,x)|^2 (1+|x|^2)^{-\left(1+\frac{m}{2}\right)}\:dx\\
&\leq & \|\phi^2 h_m\|_{L^1}^{a\theta} \| \phi^2 h_m \|_{L^{p^*}}^{a(1-\theta)}  \|  \phi h_m \|_{L^1(\R^N)}^b\\
&\leq & C \|\phi^2 h_m\|_{L^1}^{a\theta + a\frac{1-\theta}{2}} \| \nabla \phi\|_{L^2(h_m)}^{a(1-\theta)}\|  \phi h_m \|_{L^1(\R^N)}^b\\
&&+ C   \|\phi^2 h_m\|_{L^1}^a \|\phi h_m \|_{L^1(\R^N)}^b.
\end{eqnarray*}
If $N=1$, then $\theta=p^{-1}$, and straightforward computations lead to 
$$
a\theta + a \frac{1-\theta}{2}= 1 -\frac{3}{m}, \quad a (1-\theta)= \frac{2}{m}.
$$
Hence, using Young's inequality, we deduce that for all $\lambda >0$, there exists a constant $C_\lambda$ such that
\begin{eqnarray}
&&	 \int_{\R^N} |\phi(\tau,x)|^2 (1+|x|^2)^{-\left(1+\frac{m}{2}\right)}\:dx\nonumber\\
&\leq &\lambda\|\phi^2(\tau) h_m\|_{L^1} + \lambda  \| \nabla \phi(\tau)\|_{L^2(h_m)}^2 + C_\lambda \|\phi (\tau) h_m\|_{L^1(\R)}^2.\label{in:interpolation}
\end{eqnarray}

If $N\geq 2$, the calculations are similar and lead to 
$$
a\theta + a \frac{1-\theta}{2}=1- \frac{N+2}{m}, \quad a(1-\theta)= \frac{2N}{m}.
$$
Hence \eqref{in:interpolation} is also valid in this case.

Using inequality \eqref{in:tU} and choosing the parameter $\lambda$ small enough leads eventually to 
\begin{eqnarray}\label{bound_e1}
\eqref{energy1}&\leq & \frac{m-N}{16}\int_{\R^N} \left|\frac{U(\tau)}{\tU(\tau)}  \right|^2 \tU(\tau)+ \frac{1}{2}\int_{\R^N} \left|\nabla\frac{U(\tau)}{\tU(\tau)}  \right|^2 \tU(\tau)\\
&&+ C \left(\int_{\R^N} |U(\tau,x)|\:dx\right)^2\nonumber
\end{eqnarray}
for all $\tau\geq \tau_0$.

\vskip2mm
\noindent \textbf{Fourth step: control of the term \eqref{energy2}.}
\vskip1MM

\begin{rem}
We recall that \eqref{energy2}=0 if the flux $A$ is linear. Hence this step is required only in the nonlinear case.

\end{rem}

Using inequality \eqref{in:Linfty_bound}, we infer that there exists a constant $C$ such that
$$
\eqref{energy2} \leq C \int_{\R^N} |U(\tau,x)|^{1+\frac{1}{N}} \left|\nabla_x \frac{U(\tau,x)}{\tU(\tau,x)}  \right|\:dx.
$$
From now on, we treat the cases $N=1$, $N=2$, and $N\geq 3$ separately, and we set $\phi=U/\tU.$

\begin{itemize}
 \item If $N=1$, we have, for all $\tau\geq \tau_0$, 
\begin{eqnarray*}
 &&\int_{\R^N} |U(\tau,x)|^{1+\frac{1}{N}} \left|\nabla_x \frac{U(\tau,x)}{\tU(\tau,x)}  \right|dx\\&\leq&C \int_{\R} \left( \p_x \phi(\tau) h_m^{1/2} \right) \left( \phi(\tau) h_m \right)^{1/2} \left( \phi^{3/2}(\tau) h_m \right)\\
&\leq& C\| \p_x \phi(\tau)\|_{L^2(h_m)}\| U(\tau)\|_{L^1(\R)}^{1/2}\left\|  \phi^{3/2}(\tau) h_m\right\|_{L^\infty(\R)}\\
&\leq & C \| \p_x \phi(\tau)\|_{L^2(h_m)}\| U(\tau)\|_{L^1(\R)}^{1/2}\left\|  \p_x\left(\phi^{3/2}(\tau) h_m\right) \right\|_{L^1(\R)}.
\end{eqnarray*}
Moreover,
$$
 \p_x\left(\phi^{3/2} h_m\right)=\frac{3}{2}\phi^{1/2} h_m\p_x \phi - \phi^{3/2}\p_x h_m,
$$
and thus
$$
\left\|  \p_x\left(\phi^{3/2} h_m\right) \right\|_{L^1(\R)}\leq C\| \phi h_m\|_{L^1}^{1/2} \left(\|\p_x \phi\|_{L^2(h_m)} + \|\phi\|_{L^2(h_m)}\right).
$$
Eventually, we obtain, using once again \eqref{in:tU},
$$
\eqref{energy2}\leq C \| U(\tau)\|_{L^1(\R)}\left[\int_{\R} \left|\nabla\frac{U(\tau)}{\tU(\tau)}  \right|^2 \tU(\tau) + \int_{\R} \left|\frac{U(\tau)}{\tU(\tau)}  \right|^2 \tU(\tau)\right].
$$
Since 
$$
\|U(\tau)\|_{L^1}\leq \| U_{ini}\|_{L^1}\quad \forall \tau\geq 0,
$$
we infer that if $\|U_{ini}\|_1$ is sufficiently small, then
\be\label{in:NL}
\eqref{energy2}\leq \frac{m-N}{16}\int_{\R^N} \left|\frac{U(\tau)}{\tU(\tau)}  \right|^2 \tU(\tau)+ \frac{1}{2}\int_{\R^N} \left|\nabla\frac{U(\tau)}{\tU(\tau)}  \right|^2 \tU(\tau).
\ee

\item If $N=2$,  using the Sobolev embedding $W^{1,1}(\R^2)\subset L^2(\R^2),$ we obtain, for $\tau\geq \tau_0$,
\begin{eqnarray*}
 \int_{\R^N} |U(\cdot,x)|^{1+\frac{1}{N}} \left|\nabla_x \frac{U(\cdot,x)}{\tU(\cdot,x)}  \right|\:dx&\leq&C \int_{\R^2}\left(|\phi|^{3/2} h_m \right) \left(|\nabla_x \phi | h_m^{1/2}\right)\\
&\leq &C\left\| |\phi|^{3/2} h_m \right\|_{L^2(\R^2)} \| \nabla_x \phi\|_{L^2(h_m)} \\
&\leq & C \left\| \nabla\left(|\phi|^{3/2} h_m\right) \right\|_{L^1}\| \nabla_x \phi\|_{L^2(h_m)} .
\end{eqnarray*}
As is the case $N=1$, we have
$$
\left\| \nabla\left(|\phi|^{3/2} h_m\right) \right\|_{L^1(\R^2)}\leq C\| \phi h_m\|_{L^1}^{1/2} \left(\|\p_x \phi\|_{L^2(h_m)} + \|\phi\|_{L^2(h_m)}\right).
$$
Hence we are led to 
$$
\eqref{energy2}\leq C \| U(\tau)\|_{L^1(\R^2)}^{1/2}\left[\int_{\R^2} \left|\nabla\frac{U(\tau)}{\tU(\tau)}  \right|^2 \tU(\tau) + \int_{\R^2} \left|\frac{U(\tau)}{\tU(\tau)}  \right|^2 \tU(\tau)\right].
$$
Following exactly the same argument as in the case $N=1$, we deduce that if $\|U_{ini}\|_1$ is sufficiently small, then \eqref{in:NL} holds.

\item If $N\geq 3$, we have, for $\tau\geq \tau_0$, 
\begin{eqnarray*}
&& \int_{\R^N} |U(\cdot,x)|^{1+\frac{1}{N}} \left|\nabla_x \frac{U(\cdot,x)}{\tU(\cdot,x)}  \right|\:dx\\&\leq& C\int_{\R^N}|\phi|^{1+\frac{1}{N}}| \nabla_x \phi| h_m^{1+\frac{1}{N}}\\
&\leq &C\int_{\R^N}\left(|\phi| h_m\right)^{1/N} \left(|\nabla_x \phi  | h_m^{1/2}\right) \left(|\phi| h_m^{1/2}\right)\\
&\leq& C \|\phi h_m\|_{L^1(\R^N)}^{1/N} \| \nabla_x \phi\|_{L^2(h_m)}\left\| |\phi| h_m^{1/2}\right\|_{L^p(\R^N)},
\end{eqnarray*}
where the parameter $p$ is such that 
$$
\frac{1}{N} + \frac{1}{2} + \frac{1}{p}= 1,
$$
i.e. $p= (2N)/(N-2).$ Using the Sobolev embedding $H^1(\R^N)\subset L^p(\R^N)$, we have
\begin{eqnarray*}
 \left\| |\phi| h_m^{1/2}\right\|_{L^p(\R^N)}&\leq & C \left\|\nabla_x\left( |\phi| h_m^{1/2}\right)\right\|_{L^2(\R^N)}\\
&\leq & C\left(\| \nabla \phi\|_{L^2(h_m)} + \|\phi\|_{L^2(h_m)}\right).
\end{eqnarray*}
Thus, once again, we obtain
$$
\eqref{energy2}\leq C \| U(\tau)\|_{L^1(\R^N)}^{1/N}\left[\int_{\R^N} \left|\nabla\frac{U(\tau)}{\tU(\tau)}  \right|^2 \tU(\tau) + \int_{\R^2} \left|\frac{U(\tau)}{\tU(\tau)}  \right|^2 \tU(\tau)\right],
$$
and thus \eqref{in:NL} holds as long as $\| U_{ini}\|_{L^1}$ is not too large.

\end{itemize}

Gathering inequalities \eqref{energy}, \eqref{bound_e1} and \eqref{in:NL}, we infer that if $\|U_{ini}\|_1$ is sufficiently small, then for all $\tau\geq \tau_0$,
\begin{multline}
 \frac{d}{d\tau } \int_{\R^N} \left|\frac{U(\tau)}{\tU(\tau)}  \right|^2 \tU(\tau)+ \frac{m-N}{8}\int_{\R^N} \left|\frac{U(\tau)}{\tU(\tau)}  \right|^2 \tU(\tau) +  \int_{\R^N} \left|\nabla\frac{U(\tau)}{\tU(\tau)}  \right|^2 \tU(\tau)\leq\\
\leq C\left(\int_{\R^N} |U(\tau,x)|\:dx\right)^2\leq C \|U_{ini}\|_{L^1}^2.
\end{multline}
\vskip2mm

\noindent \textbf{Fifth step: Conclusion.}
\vskip1mm

Let $C_1:=(m-N)/4$, $C_2:= C \|U_{ini}\|_{1}^2$. Using a Gronwall type argument, we deduce that for all $\tau\geq \tau_0,$ we have
\begin{eqnarray*}
 &&\int_{\R^N} \left|\frac{U(\tau)}{\tU(\tau)}  \right|^2 \tU(\tau) + \int_{\tau_0}^\tau e^{-C_1(\tau-s)}\left(\int_{\R^N} \left|\nabla\frac{U(s)}{\tU(s)}  \right|^2 \tU(s)\right)\:ds\\
&\leq& e^{-C_1(\tau-\tau_0)}\int_{\R^N} \left|\frac{U(\tau_0)}{\tU(\tau_0)}  \right|^2 \tU(\tau_0) + \frac{C_2}{C_1}\\
&\leq & C \int_{\R^N} |U(\tau_0, x)|^2 (1+ |x|^2)^{m/2} + \frac{C_2}{C_1}.
\end{eqnarray*}
Using \eqref{in:tU}, we infer
$$\begin{aligned}
  \sup_{\tau\geq \tau_0} \int_{\R^N} |U(\tau, x)|^2 (1+|x|^2)^{m/2}\leq C \int_{\R^N} |U(\tau_0, x)|^2 (1+ |x|^2)^{m/2} +  C \frac{C_2}{C_1},\\
	\sup_{\tau\geq \tau_0} \int_{\tau}^{\tau+1}\int_{\R^N} \left|\nabla\frac{U(s)}{\tU(s)}  \right|^2 h_m\:ds\leq C \int_{\R^N} |U(\tau_0, x)|^2 (1+ |x|^2)^{m/2} +  C \frac{C_2}{C_1}.
  \end{aligned}
$$
Hence $U\in L^\infty([0,\infty), L^2((1+|x|^2)^{m/2})).$
Since $f_0$ is bounded away from zero, the $L^2$ bound on $V$ follows. 

Concerning the bound on $\nabla_x V$, notice that
$$
V(\tau,x)=\frac{U(\tau,x)}{\tU(\tau,x)}\left(h_m(x) + e^{-t} \frac{f_1(z)}{f_0(z)}\cdot \nabla_x h_m(x)\right), 
$$
and thus
$$
\left| \nabla_x V(\tau,x)\right|\leq C h_m(x)\left(\left|\nabla_x \frac{U(\tau,x)}{\tU(\tau,x)}\right| + \left| \frac{U(\tau,x)}{\tU(\tau,x)}\right| \right).
$$
Consequently, for all $\tau\geq 0$
$$
\int_{\R^N}|\nabla_x V(\tau,x)|^2 (1+|x|^2)^{m/2}\:dx\leq C \left\| \frac{U(\tau)}{\tU(\tau)}  \right\|_{H^1(h_m)}^2,
$$
which leads to the bound on $\nabla_x V$. (Notice that we even retrieve that $\nabla_x V \in L^2_\text{loc}([0,\infty), L^2((1+|x|^2)^{m/2}))$). 

Hence Proposition \ref{prop:compactness} is proved.

\vskip2mm

Let us now conclude this section by explaining how the bound \eqref{hyp:moments} on the moments of order four follows from \eqref{in:L2control}. Let $U_{ini}\in L^2(h_m^{-1})$, with $m>2(N+2)$ sufficiently large. Then we have proved that $U\in L^\infty([0,\infty), L^2(h_m^{-1}))$. Now, for all $\tau\geq 0$, using a simple H\"older inequality, we infer that
$$
\int_{\R^N} |U(\tau,x)| |x|^4\:dx\leq \|U(\tau)\|_{L^2(h_m^{-1})} \left(\int_{\R^N} |x|^8 (1+|x|^2)^{-m/2}\:dx\right)^{1/2}.
$$
Hence, if $m>N+8$, we deduce that $U\in L^\infty([0,\infty), L^1(|x|^4))$; going back to the original variables, this entails that \eqref{hyp:moments} is satisfied. Thus the convergence result \eqref{in:BDK} holds if the flux $A$ is linear, and Proposition \ref{prop:linear} is proved.

\section{Long-time behaviour}

\label{sec:conclusion}

This section is devoted to the rest of the proof of Theorem \ref{thm:conv}. The idea is to use the $L^1$ compactness proved in the previous section (see Proposition \ref{prop:compactness}) together with techniques from dynamical systems theory. This type of proof was  initiated by S. Osher and J. Ralston in \cite{OR}, in which the authors proved the $L^1$ stability of travelling waves for a quasilinear parabolic equation. Their arguments were then adapted successfully to various kinds of problems in the context of scalar conservation laws (see for instance the review in \cite{SerreHandbook}).

In the present study, our scheme of proof is in fact closely related to the one of M. Escobedo and E. Zuazua in \cite{EZ}; indeed, the idea is to apply the dynamical systems tools to the rescaled parabolic system \eqref{eq:rescaled} rather than the original conservation law \eqref{LCS}. The main difference with \cite{EZ} lies in the presence of highly oscillating coefficients in \eqref{eq:rescaled}; thus it is necessary to work simultaneously with the homogenized equation \eqref{eq:h_Ngrand}-\eqref{eq:h_N=1} and with the oscillating one.

Let us now introduce some notation and definitions. First, we denote by $S_\tau$ ($\tau\geq 0$) the semi-group associated with the homogenized equation, that is equation \eqref{eq:h_N=1} if $N=1$, and equation \eqref{eq:h_Ngrand} if $N\geq 2$. According to Proposition \ref{prop:ex_uni_h}, the semi-group $S_\tau$ is well-defined in $L^\infty(\R^N)\cap L^2(K)$; additionally, the $L^1$ contraction property holds, namely
$$
\| S_\tau F_1- S_\tau F_2\|_{L^1(\R^N)}\leq \| F_1- F_2\|_{L^1(\R^N)}\quad \forall \tau\geq 0,\ \forall F_1, F_2\in L^\infty(\R^N)\cap L^2(K).
$$
Hence $S_\tau$ can be extended on $L^1(\R^N)$.

We also define the $\omega$-limit set associated with a given function $U_{ini}\in L^1(\R^N)$: recalling the definition of the function $V$ (see \eqref{def:V}), we set
$$
\Omega[U_{ini}]:=\left\{ \bV \in L^1(\R^N),\quad \exists \tau_n\to\infty,\ V(\tau_n)\to \bV\quad\text{in } L^1(\R^N)\right\},
$$
where the function $U$ in \eqref{def:V} is the unique solution of \eqref{eq:rescaled} with initial data $U_{ini}$. When there is no ambiguity, we will simply write $\Omega$ instead of $\Omega[U_{ini}]$.

Notice that $V(\tau_n)$ converges towards $\bV$ in $L^1$ if and only if
$$
\lim_{n\to\infty}\int_{\R^N} \left| U(\tau_n, x) - f_0\left(e^{\tau_n} x + c \frac{e^{2\tau_n} -1}{2}\right) \bV(x)\right|\:dx=0.
$$
This equivalence will be used repeatedly throughout the section.

The organisation of this section is the following: we first introduce a ``quasi-Lyapunov function'' for the semi-group associated with equation \eqref{eq:rescaled}. We then prove that Proposition \ref{prop:conv_rescaled} holds when the initial data $U_{ini}$ has a sufficiently small $L^1$ norm. Eventually, we prove Proposition \ref{prop:conv_rescaled} in the general case.

\subsection{A quasi-Lyapunov function}

Let us introduce the notion of quasi-Lyapunov function:

\begin{defi}
 Let $\mathcal X$ be a Banach space, and let $\mathcal H:[0,\infty)\times \mathcal X\to \R$. Let $\{u(t)\}_{t\geq 0}$ be a trajectory in $\mathcal X$. We say that $\mathcal H$ is a quasi-Lyapunov function for the trajectory $u$ if the following properties hold:

\begin{enumerate}[(i)]
 \item The family $\mathcal H(t,u(t))$ ($t\geq 0$) is bounded in $\R$;
\item There exists  a function $\psi:[0,\infty)\to[0,\infty)$ such that $\lim_{t\to\infty}\psi(t)=0$ and
$$
\forall t\geq 0,\quad \sup_{s\geq t} (\mathcal H(s, u(s))-\mathcal H(t, u(t)))\leq \psi(t).
$$

\end{enumerate}

\label{def:q-Lyapunov}

\end{defi}

We then have the following result: 
\begin{lem}
 Let $\mathcal X$ be a Banach space, and let $\{u(t)\}_{t\geq 0}$ be a trajectory in $\mathcal X$. Let $\mathcal H:[0,\infty)\times \mathcal X\to \R$ be a quasi-Lyapunov function for the trajectory $u$. Then $\mathcal H(t,u(t))$ has a finite limit as $t\to\infty.$

\label{lem:Lyapunov}
\end{lem}

\begin{proof}
First, since $\mathcal H(t,u(t))$ is bounded for $t\in[0,\infty)$, the quantities
$$
\underline{H}:=\liminf_{t\to\infty} \mathcal H(t,u(t)),\quad \overline{H}:=\limsup_{t\to\infty} \mathcal H(t,u(t))
$$ 
are well-defined and belong to $\R$, with $\underline{H}\leq \overline{H}$.

Let $\e>0$ arbitrary. There exists $t_\e>0$ such that
$$
\psi(t)\leq \e\quad\forall t\geq t_\e.
$$
By definition of $\underline{H}$, there exists $s_\e\geq t_\e$ such that
$$
\left|\mathcal H(s_\e,u(s_\e))- \underline{H}\right|\leq \e.
$$
Since $\mathcal H$ is a quasi-Lyapunov function, for all $s\geq s_\e$, we have
\begin{eqnarray*}
\mathcal H(s,u(s))&\leq &  \mathcal H(s_\e,u(s_\e)) + \psi(s_\e)\\
&\leq & \underline{H} + 2\e.
\end{eqnarray*}
Hence
$$
\overline{H}\leq \underline{H} + 2\e\quad\forall \e>0,
$$
and $\overline{H}=\underline{H} $. Thus the quantity $\mathcal H(t,u(t))$ has a finite limit as $t\to\infty$.

\end{proof}

We now apply this notion to the present context:

\begin{lem}
Let $M\in\R$ be arbitrary, and let $U_{ini}\in L^1(\R^N)$. For $\tau\geq 0$ and $u\in L^1(\R^N)$, define the function $\mathcal H$ by
$$
\mathcal H(\tau,u):=\int_{\R^N} \left|u(x)- \Uapp[F_M](\tau,x;e^\tau)\right|\:dx,
$$
where the function $\Uapp$ was introduced in Definition \ref{def:Uapp}.

Let $U\in\mathcal C(|0,\infty), L^1(\R^N))$ be  the solution of \eqref{eq:rescaled} with initial data $U_{\tau=0}=U_{ini}.$

Then $\mathcal H$ is a quasi-Lyapunov function for the trajectory $\{U(\tau))\}_{\tau\geq 0}$ in $L^1(\R^N)$. As a consequence, the function
$$
\tau\in[0,\infty)\mapsto \int_{\R^N} |U(\tau,x) - f_0(z) F_M(x)|\:dx,\quad \text{with }z=e^\tau x + c\frac{e^{2\tau} -1}{2},
$$
converges as $\tau\to \infty.$

\end{lem}

\begin{proof}

This property is an easy consequence of the first point in Lemma \ref{lem:approx_sol}; indeed, according to Lemma \ref{lem:approx_sol}, there exists a constant $C$, depending only on $N$ and $M$, such that
$$
\frac{d}{d\tau} \mathcal H(\tau,U(\tau))= \frac{d}{d\tau}\| U(\tau) - \Uapp[F_M](\tau)\|_{L^1(\R^N)}\leq C e^{-\tau}.
$$
Consequently, for all $\tau'\geq \tau\geq 0,$ we have
$$
\mathcal H(\tau', U(\tau')) -\mathcal H(\tau, U(\tau))\leq C(e^{-\tau} - e^{-\tau'})\leq C e^{-\tau}.
$$
Thus property (ii) of Definition \ref{def:q-Lyapunov} is satisfied. Additionally, notice that
\begin{eqnarray*}
0\leq \mathcal H(\tau,U(\tau))&\leq& \| U(\tau)\|_{L^1(\R^N)} + \|f_0\|_{L^\infty(\T^N)}\| F_M\|_{L^1(\R^N)} \\&+& Ce^{-\tau} \left(\|\nabla F_M\|_{L^1(\R^N)}+\|F_M\|_{L^2(\R^N)}^2 \right)\\
&+&C e^{-2\tau}\left(\|F_M\|_{W^{2,1}(\R^N)} + \| F_M\|_{L^3(\R^N)}^3 + \|\nabla F_M\|_{L^2(\R^N)}^2\right)\\
&\leq & \|U_{ini}\|_{L^1(\R^N)} + C.
\end{eqnarray*}
Whence $\mathcal H(\tau,U(\tau))$ is bounded for $\tau\in[0,\infty)$. Consequently $\mathcal H$ is a quasi-Lyapunov function for the trajectory $U(\tau)$. According to Lemma \ref{lem:Lyapunov}, $\mathcal H(\tau,U(\tau))$ admits a finite limit as $\tau\to\infty.$ Furthermore, we have
$$
\left\| U(\tau) - f_0(z) F_M\right\|_{L^1}= \left\| U(\tau) - \Uapp[F_M](\tau)  + e^{-\tau} U_1 + e^{-2\tau} U_2\right\|_{L^1},
$$
where $U_1$ and $U_2$ are defined by \eqref{def:U_1_Ngrand}-\eqref{def:U_1_N=1} and \eqref{eq:U2} respectively. Hence for all $\tau\geq 0$, there holds
$$
\mathcal H(\tau,U(\tau)) - C e^{-\tau} \leq \left\| U(\tau) - f_0(z) F_M\right\|_{L^1} \leq \mathcal H(\tau,U(\tau)) + C e^{-\tau},
$$
where the constant $C$ depends only on $W^{s,p}$ bounds on $F_M$. Thus the function
$$
\tau\mapsto \left\| U(\tau) - f_0(z) F_M\right\|_{L^1(\R^N)} 
$$
converges as $\tau\to\infty,$ and 
$$
\lim_{\tau\to\infty} \left\| U(\tau) - f_0(z) F_M\right\|_{L^1} = \lim_{\tau\to\infty} \mathcal H(\tau,U(\tau)).
$$

\end{proof}

\begin{defi}
Let $U_{ini}\in L^1(\R^N) $ be arbitrary, and let $M:=\int_{\R^N} U_{ini}$. Let $U$ be the solution of \eqref{eq:rescaled} with initial data $U_{|t=0}=U_{ini}$.

We define the number $\ell(U_{ini})$ by
$$
\ell(U_{ini}):=\lim_{\tau\to\infty}\int_{\R^N} |U(\tau,x) - f_0(z) F_M(x)|\:dx,\quad \text{with }z=e^\tau x + c\frac{e^{2\tau} -1}{2}.
$$

\end{defi}
Notice that Proposition \ref{prop:conv_rescaled} is equivalent to
$$
\ell(U_{ini})=0\quad\forall U_{ini}\in L^1(\R^N).
$$

Classically, we now derive a continuity property for the function $\ell$:
\begin{lem}
The function
$$
U\in L^1(\R^N)\mapsto \ell(U)\in\R
$$ 
is Lipschitz continuous.

\label{lem:ell_cont}
\end{lem}

\begin{proof}
Let $U_{ini}^{(1)},U_{ini}^{(2)}\in L^1(\R^N) $, and let $M^{(i)}= \int_{\R^N} U_{ini}^{(i)}$ for $i=1,2$. We denote by $U^{(i)}\in\mathcal C([0,\infty), L^1(\R^N))$ the solution of \eqref{eq:rescaled} with initial data $U_{ini}^{(i)}$. Then for all $\tau\geq 0$, the $L^1$ contraction principle ensures that
$$
\left\| U^{(1)}(\tau) - U^{(2)}(\tau)\right\|_{L^1(\R^N)}\leq \left\| U_{ini}^{(1)} - U_{ini}^{(2)}\right\|_{L^1(\R^N)}.
$$
Hence, for all $\tau\geq 0$, we have
\begin{eqnarray*}
&&\left| \int_{\R^N} \left|U^{(1)}(\tau,x) - f_0(z) F_{M^{(1)}}(x)  \right|dx -  \int_{\R^N} \left|U^{(2)}(\tau,x) - f_0(z) F_{M^{(2)}}(x)  \right|dx  \right|\\
&\leq & \left\| U^{(1)}(\tau) - U^{(2)}(\tau)\right\|_{L^1(\R^N)} + \left\| f_0\right\|_{L^\infty(\T^N)}\left\| F_{M^{(1)}}-F_{M^{(2)}}\right\|_{L^1(\R^N)}.
\end{eqnarray*}
According to Lemma A.2 in the Appendix, 
$$
\left\| F_{M^{(1)}}-F_{M^{(2)}}\right\|_{L^1(\R^N)}=\left| M^{(1)}-  M^{(2)} \right|\leq  \left\| U_{ini}^{(1)} - U_{ini}^{(2)}\right\|_{L^1(\R^N)}.
$$
Eventually, we obtain, for all $\tau\geq 0$,
\begin{multline*}
\left| \int_{\R^N} \left|U^{(1)}(\tau,x) - f_0(z) F_{M^{(1)}}(x)  \right|dx -  \int_{\R^N} \left|U^{(2)}(\tau,x) - f_0(z) F_{M^{(2)}}(x)  \right|dx  \right|\leq \\\leq\left(1+\left\| f_0\right\|_{L^\infty(\T^N)}\right)\left\| U_{ini}^{(1)} - U_{ini}^{(2)}\right\|_{L^1(\R^N)},
\end{multline*}
and thus, passing to the limit,
$$
\left| \ell\left( U_{ini}^{(1)} \right) - \ell\left( U_{ini}^{(2)} \right) \right|\leq \left(1+\left\| f_0\right\|_{L^\infty(\T^N)}\right)\left\| U_{ini}^{(1)} - U_{ini}^{(2)}\right\|_{L^1(\R^N)}.
$$
Hence $\ell$ is a Lipschitz continuous function.

\end{proof}

\subsection{Analysis of the $\omega$-limit set}

\begin{prop} Let $U_{ini}\in L^1(\R^N)$, and set 
 $$
M:=\int_{\R^N} U_{ini}.
$$

Assume that the $\omega$-limit set $\Omega$ associated with $U_{ini}$ is non-empty. Then the following properties hold:
\begin{enumerate}[(i)]
\item For all $\bar V\in\Omega$, 
$$
\int_{\R^N}\bV=M;
$$

\item $S_\tau \Omega\subset \Omega$ for all $\tau\geq 0$;

\item For all $\bV\in\Omega$, we have 
$$
\left\| \bV -F_M \right\|_{L^1(\R^N)}=\ell(U_{ini}).
$$

\end{enumerate}

\label{prop:dynamical}
\end{prop}

\begin{proof}
Throughout the proof, we denote by $U$ the unique solution of equation \eqref{eq:rescaled} with initial data $U_{ini}$.

Property (i) is  quite straightforward: indeed, conservation of mass for the equation \eqref{eq:rescaled} implies that
$$
\int_{\R^N} U(\tau)=M\quad\forall \tau\geq 0.
$$
If $\bV\in\Omega$, then there exists a sequence $(\tau_n)_{n\geq 0}$ such that
$$
\lim_{n\to\infty}\tau_n=\infty\text{ and }\lim_{n\to\infty}\int_{\R^N} \left| U(\tau_n, x) - f_0\left(z_n\right) \bV(x)\right|\:dx=0,
$$
where $z_n=e^{\tau_n} x + c \frac{e^{2\tau_n} -1}{2}$. According to a result of G. Allaire (see \cite{Allaire}), 
$$
\lim_{n\to\infty}\int_{\R^N}f_0\left(z_n\right) \bV(x)\:dx=\mean{f_0}\int_{\R^N} \bV=\int_{\R^N} \bV;
$$
gathering the three equalities, we obtain property (i).

\vskip1mm

We now address the proof of property (ii), which relies on the second point in Lemma \ref{lem:approx_sol}; let $\bar V\in\Omega$ be arbitrary, and for all $\e>0$, let $\bV_\e\in L^2(K)\cap L^\infty(\R^N)$ such that 
$$
\| \bV_\e -\bV\|_{L^1(\R^N)}\leq \e.
$$
Let $(\tau_n)_{n\geq 0}$ be a sequence of positive numbers such that $\tau_n\to\infty$ and 
$$
\int_{\R^N} \left| U(\tau_n,x) - f_0(z_n) \bV(x)\right|\:dx\to 0,
$$
where $z_n=e^{\tau_n} x + c \frac{e^{2\tau_n}-1}{2}$.

Let $\rho\in\mathcal C^\infty_0(\R^N)$ be a mollyfing kernel; for $\delta>0$, set $\rho_\delta:=\delta^{-N}\rho(\cdot/\delta)$, and define the function $\unde$ by
$$
\unde(\tau,x):=\Uapp\left[ (S_\tau \bV_\e)\ast_x\rho_\delta\right](\tau_n+ \tau,x; e^{\tau_n+\tau}).
$$
Then Lemma \ref{lem:approx_sol} ensures that $\unde$ satisfies equation \eqref{eq:rescaled} with an error term, the latter being   bounded for all $T>0$ in $L^\infty([0,T], L^1(\R^N))$ by 
$$
\omega_{T,\e}(\delta) + C_{T,\e,\delta} e^{-\tau_n}
$$
where $\omega_{T,\e}:[0,\infty)\to[0,\infty)$ is such that $\lim_{0^+} \omega_{T,\e}=0,$ and where the constant $C_{T,\e,\delta}$ depends only on $\e,\delta,N$ and $T$.

Using the $L^1$ contraction principle for scalar conservation laws, we infer that for all $T>0$, and for all $\tau\in[0,T]$, 
\begin{eqnarray*}
&&\int_{\R^N}\left|U(\tau_n+ \tau, x) - U_n^{\delta,\e}(\tau,x)\right|\:dx\\
&\leq & \omega_{T,\e}(\delta) + C_{T,\e,\delta} e^{-\tau_n} + \int_{\R^N}\left|U(\tau_n, x) - U_{n|\tau=0}^{\delta,\e}(x)\right|\:dx\\
&\leq & \omega_{T,\e}(\delta) + C_{T,\e,\delta} e^{-\tau_n} + \int_{\R^N} \left| U(\tau_n,x) - f_0(z_n) \bV(x)\right|\:dx\\
&&+ \int_{\R^N} \left|  U_{n|\tau=0}^{\delta,\e}(x) - f_0(z_n) \bV(x)\right|\:dx.
\end{eqnarray*}
Now, according to Definition \ref{def:Uapp},
\begin{eqnarray*}
 U_n^{\delta,\e}(\tau,x)&=& (S_\tau \bV_\e)\ast_x \rho_\delta(x) f_0\left(e^{\tau_n+\tau} x + c \frac{e^{2(\tau_n + \tau)}-1}{2}\right) \\&+& e^{-(\tau_n+\tau)} U_1\left( \tau,x, e^{\tau_n+\tau} x+c \frac{e^{2(\tau_n + \tau)}-1}{2}\right) \\&+&   e^{-2(\tau_n+\tau)}U_2\left( \tau,x, e^{\tau_n+\tau} x + c \frac{e^{2(\tau_n + \tau)}-1}{2}\right).
\end{eqnarray*}
Hence for all $\tau\in[0,T]$, we have
\begin{eqnarray*}
&&\int_{\R^N}\left| U_n^{\delta,\e}(\tau,x) - S_\tau \bV f_0 \left(e^{\tau_n+\tau} x + c \frac{e^{2(\tau_n + \tau)}-1}{2}\right) \right|dx\\&\leq & \left\| f_0  \right\|_\infty \sup_{\tau\in[0,T]} \| S_\tau\bV - (S_\tau \bV_\e)\ast_x \rho_\delta\|_{L^1(\R^N)}\\
&&+ e^{-\tau_n} \left( \| U_1\|_{L^\infty([0,T]\times \T^N_z, L^1(\R^N_x))} + \| U_2\|_{L^\infty([0,T]\times \T^N_z, L^1(\R^N_x))}  \right)\\
&\leq & \left\| f_0  \right\|_\infty \left[\| \bV-\bV_\e\|_{L^1} + \sup_{\tau\in[0,T]} \| S_\tau\bV_\e - (S_\tau \bV_\e)\ast_x \rho_\delta\|_{L^1(\R^N)}\right]\\
&&+e^{-\tau_n} \left( \| U_1\|_{L^\infty([0,T]\times \T^N_z, L^1(\R^N_x))} + \| U_2\|_{L^\infty([0,T]\times \T^N_z, L^1(\R^N_x))}  \right)\\
&\leq & C\e + \omega_{T,\e}(\delta) + C_{T,\e,\delta} e^{-\tau_n}.
\end{eqnarray*}
Gathering the two inequalities, we deduce that
for all $\tau\in [0,T]$, for all $n,\delta,\e$,
\begin{eqnarray*}
&&\int_{\R^N}\left|U(\tau_n+ \tau, x) - S_\tau \bV (x)f_0 \left(e^{\tau_n+ \tau} x + c \frac{e^{2(\tau_n+\tau)}-1}{2}\right)\right|\:dx\\
&\leq & \omega_{T,\e}(\delta) + C_{T,\e,\delta} e^{-\tau_n} + \int_{\R^N} \left| U(\tau_n,x) - f_0(z_n) \bV(x)\right|\:dx + C\e.
\end{eqnarray*}In the right-hand side of the above inequality, we first choose $\e$ sufficiently small, then $\delta$ so that $\omega_{T,\e}(\delta)$ is sufficiently small, and eventually $n$ large enough so that the two remaining terms are small as well; hence
$$
\lim_{n\to\infty} \inf_{\e>0,\delta>0}\left( \omega_{T,\e}(\delta) + C_{T,\e,\delta} e^{-\tau_n} + \int_{\R^N} \left| U(\tau_n,x) - f_0(z_n) \bV(x)\right|\:dx + C\e\right)=0.
$$
Thus we have proved that for all $T>0$, 
$$
\lim_{n\to\infty}\sup_{\tau\in[0,T]}\int_{\R^N}\left|U(\tau_n+ \tau, x) - (S_\tau \bV) f_0 \left(e^{\tau_n+ \tau} x + c \frac{e^{2(\tau_n+\tau)}-1}{2}\right)\right|\:dx=0.
$$
The above convergence entails immediately that $S_\tau \bV\in \Omega$ for all $\tau\in[0,T]$. Since $T>0$ was arbitrary, property (ii) is proved.

\vskip1mm

There remains to prove property (iii), which is a variant of the LaSalle invariance principle; let $\bar V\in\Omega$ be arbitrary, and let $\tau_n$ be a sequence of positive numbers such that $\lim_{n\to\infty} \tau_n=+\infty$ and
$$
\lim_{n\to\infty}\int_{\R^N}\left| U(\tau_n,x) - f_0(z_n) \bV(x)\right|\:dx=0, 
$$
where $z_n= e^{\tau_n} x + c \frac{e^{2\tau_n}-1}{2}$. According to a result of G. Allaire (see \cite{Allaire}), we have, since $\mean{f_0}=1$ and $f_0\in\mathcal C(\T^N)$,
\begin{eqnarray*}
\left\| \bV -F_M\right\|_{L^1(\R^N)} &=&\lim_{n\to\infty} \int_{\R^N} f_0(z_n) |\bV(x) -F_M(x)|\:dx\\
&=&\lim_{n\to\infty} \int_{\R^N} f_0(z_n) \left|\bV(x) - \frac{U(\tau_n, x)}{f_0(z_n)} +\frac{U(\tau_n, x)}{f_0(z_n)}   -F_M(x)\right|dx\\
&=&\lim_{n\to\infty} \int_{\R^N} f_0(z_n) \left|\frac{U(\tau_n, x)}{f_0(z_n)}   -F_M(x)\right|\:dx\\
&=&\ell(U_{ini}).
\end{eqnarray*}
Consequently,
$$
\| \bV -F_M\|_{L^1(\R^N} =\ell(U_{ini})\quad\forall \bV\in\Omega.
$$
\end{proof}

\begin{corol}
  Let $U_{ini}\in L^1(\R^N)$, and set 
 $$
M:=\int_{\R^N} U_{ini}.
$$

Assume that the $\omega$-limit set $\Omega[U_{ini}]$ is non-empty. Then $\ell (U_{ini})=0$, and thus the result of Proposition \ref{prop:conv_rescaled} holds.

\label{cor:omega/conv}
\end{corol}

\begin{proof}
Let $\bV\in \Omega$ be arbitrary. Then 
$$
\lim_{\tau\to\infty}\| S_\tau\bV -F_M\|_{L^1(\R^N)}=0;
$$
this property is stated in Proposition \ref{prop:ex_uni_h} in the case when $\bar V\in L^\infty(\R^N)\cap L^2(K)$, but can be in fact easily generalized to an arbitrary function $\bV\in L^1$ by using the contractivity of the semi-group $S_\tau$: indeed, let $\e>0$, and let $\bV_\e\in  L^\infty(\R^N)\cap L^2(K)$ such that $\int \bV_\e=\int \bV=M $, and 
$$
\|\bV_\e -\bV \|_{L^1(\R^N)}\leq \e.
$$
Then for all $\tau\geq 0$,
$$
\| S_\tau \bV - F_M\|_1\leq \| S_\tau \bV-S_\tau \bV_\e\|_1 + \| S_\tau \bV_\e - F_M\|_1\leq \| \bV-\bV_\e\|_1+ \| S_\tau \bV_\e - F_M\|_1.
$$
Hence, using Proposition \ref{prop:ex_uni_h}, we infer that
$$
\limsup_{\tau\to\infty} \| S_\tau \bV - F_M\|_1\leq\e\quad\forall \e>0,
$$
and thus $\| S_\tau\bV -F_M\|_1$ vanishes as $\tau\to\infty$.

On the other hand, property (ii) in Proposition \ref{prop:dynamical} ensures that $S_\tau\bV\in\Omega$ for all $\tau\geq 0$, and thus, using (iii),
$$
\|S_\tau\bV -F_M\|_1=\ell(U_{ini})\quad\forall \tau\geq 0.
$$
Consequently, $\ell(U_{ini})=0.$ Going back to the definition of $\ell(U_{ini})$, we deduce that
$$
\lim_{\tau\to\infty}\int_{\R^N} \left| U(\tau,x) - f_0\left(e^\tau x + c \frac{e^{2\tau}-1}{2}  \right)F_M(x) \right|\:dx=0.
$$
\end{proof}

Thus the proof of Proposition \ref{prop:conv_rescaled} is complete provided we are able to show that the set $\Omega[U_{ini}]$ is non-empty for a sufficiently large class of functions $U_{ini}\in L^1(\R^N)$. In the case when $\|U_{ini}\|_1$ is small, this result follows from Proposition \ref{prop:compactness} and from a contraction principle. The proof in the general case is more involved, and in fact, an analysis similar to the one performed in Section \ref{sec:bounds} has to be conducted once more.

\subsection{Proof of Proposition \ref{prop:conv_rescaled} when $\| U_{ini}\|_1$ is small}

We now complete the proof of Theorem \ref{thm:conv} when $\|U_{ini}\|_1$ is small. Let $U_{ini}\in L^1(\R^N)$. Assume that   $U_{ini}$ satisfies the following assumptions
\begin{eqnarray}
\exists m>2(N+1),&&\|U_{ini}\|_{L^1(\R^N)}\leq C_m,\label{hyp:U_{ini}_petit}\\
\text{and}&&U_{ini}\in L^2((1+|x|^2)^{m/2})\cap L^\infty(\R^N),\label{hyp:U_{ini}_int}
\end{eqnarray}
where the constant $C_m$ was introduced in Proposition \ref{prop:compactness}. Then according to Proposition \ref{prop:compactness}, the $\omega$-limit set $\Omega[U_{ini}]$ is non-empty, and consequently Proposition \ref{prop:conv_rescaled} is true (see Corollary \ref{cor:omega/conv}).

Let us now prove that Proposition \ref{prop:conv_rescaled} holds when $U_{ini}$ merely satisfies \eqref{hyp:U_{ini}_petit}: this fact is a direct consequence of the density of $L^2((1+|x|^2)^{m/2})\cap L^\infty(\R^N)$ in $L^1(\R^N)$, together with the continuity of $\ell$. Indeed, for all $\e>0$, let $U_{ini}^\e\in L^2((1+|x|^2)^{m/2})\cap L^\infty(\R^N)$ such that
$$
\| U_{ini} - U_{ini}^\e\|_{L^1(\R^N)} \leq \e,\quad\|U_{ini}^\e\|_{L^1(\R^N)}\leq C_{m}.
$$
Then $\ell(U_{ini}^\e)=0.$ Since $\ell$ is Lipschitz continuous (see Lemma \ref{lem:ell_cont}), there exists a constant $C$ such that
$$
\ell(U_{ini})= |\ell(U_{ini})-\ell(U_{ini}^\e)|\leq C \| U_{ini} - U_{ini}^\e\|_{L^1(\R^N)}\leq C\e.
$$
Since the above inequality holds for all $\e>0$, we deduce that $\ell(U_{ini})=0.$ Recalling the definition of $\ell$, we infer that Proposition \ref{prop:conv_rescaled} holds for all initial data $U_{ini}\in L^1(\R^N)$ satisfying \eqref{hyp:U_{ini}_petit}.

\subsection{Proof of Proposition \ref{prop:conv_rescaled} in the general case}

The case when $\|U_{ini}\|_1$ is large follows from the  following Lemma:

\begin{lem}
There exists a constant $C_0$, depending only on $N$ and on the flux $A$, such that for all $U_{ini}\in L^1(\R^N)$,
$$
\ell(U_{ini})\leq C_0 \Rightarrow \ell(U_{ini})=0.
$$ 
\label{lem:lU0}
\end{lem}

Before proving the above Lemma, let us explain why the result of Proposition \ref{prop:conv_rescaled} follows. If $\ell(U_{ini})\leq C_0$ for all $U_{ini}\in L^1(\R^N)$, then the above Lemma states that $\ell$ is identically zero on $L^1(\R^N)$, and thus Proposition \ref{prop:conv_rescaled} is true.  Thus we assume by contradiction that there exists $U_{ini}\in L^1(\R^N)$ such that $\ell(U_{ini})>C_0$. Consider the function
$$
\phi:t\in[0,1]\mapsto \ell(t U_{ini}).
$$
We have proved in the previous paragraph that $\phi$ vanishes in a neighbourhood of zero. Moreover, $\phi$ is a continuous function according to Lemma \ref{lem:ell_cont}.
Now, it is obvious that $\phi(0)=0$, and $\phi(1)=\ell(U_{ini})>C_0.$ Hence there exists $t_0\in(0,1)$ such that
$$
\phi(t_0)=\frac{C_0}{2}
$$
But according to Lemma \ref{lem:lU0}, $\phi(t_0)=\ell(t_0 U_{ini})=0$, which is impossible. Thus $\ell(U_{ini})=0$ for all $U_{ini}\in L^1(\R^N)$.

There remains to prove Lemma \ref{lem:lU0}. According to Corollary \ref{cor:omega/conv} and using by now standard arguments, we only have to prove that there exists a set $\mathcal A\subset L^1(\R^N)$, which is dense in $L^1(\R^N)$, and such that
\be\label{implication}
\exists C>0,\ \forall U_{ini}\in L^1(\R^N)\cap \mathcal A,\quad \ell(U_{ini})\leq C \Rightarrow \Omega[U_{ini}]\neq \emptyset.
\ee
In the following, we will take $\mathcal A=L^2((1+|x|^2)^{m/2})$, for some $m>0$ sufficiently large.

The scheme of proof of the implication \eqref{implication} is very similar to the one of Proposition \ref{prop:compactness}; indeed, we have to prove that if $\ell(U_{ini})$ is small enough, then there exists a sequence $(\tau_n)$ of positive numbers, with $\lim_{n\to\infty}\tau_n=+\infty$, such that $(V(\tau_n,\cdot))_{n\to\infty}$ is a compact sequence in $L^1(\R^N)$. Notice that this is obviously equivalent to the compactness of the sequence $V(\tau_n,\cdot) - F_M$, whose $L^1$ norm is of the order of $\ell(U_{ini})$ as $n\to\infty$. Thus our strategy is the following: rather than using directly the equation on $U$, we consider the equation on the function $U-\Uapp[F_M]$. We prove that for an appropriate function $\tilde U$, an inequality of the type \eqref{energy} holds, with $U$ replaced by $U-\Uapp[F_M]$. Then, all the occurrences of $\| U(\tau)\|_1$ in the proof of Proposition \ref{prop:compactness} are replaced by $\|(U-\Uapp[F_M])(\tau)\|_1$, which converges towards $\ell(U_{ini})$ as $\tau\to\infty$. Thus the same arguments which led us to compactness in the case when $\|U_{ini}\|_1$ is small show that compactness holds, provided $\ell(U_{ini})$ is small enough.

Let us now retrace the main lines of the proof: first, consider a function $U_{ini}\in L^1(\R^N)$ such that $U_{ini}\in L^2((1+|x|^2)^{m/2}) $ for some sufficiently large $m$ (to be chosen later). Set $M=\int_{\R^N} U_{ini}$ and
$$
W(\tau,x)= U(\tau,x) - \Uapp [F_M](\tau,x; e^{\tau}).
$$
In the rest of the proof, for the sake of brevity, we will write $\Uapp(\tau,x)$ as a short-hand for $\Uapp [F_M](\tau,x; e^{\tau})$.
Then the following properties hold
$$
\begin{aligned}
W\in L^\infty_\text{loc}([0,\infty),   L^2((1+|x|^2)^{m/2}))\cap L^2_\text{loc}([0,\infty),   H^1((1+|x|^2)^{m/2})),\\
\exists C>0,\ \forall\tau\geq 0,\quad \|W(\tau,\cdot)\|_{L^\infty(\R^N)}\leq C e^{N\tau},\\
\lim_{\tau\to\infty}\|W(\tau)\|_{L^1(\R^N)}=\ell(U_{ini}).
\end{aligned}
$$
Moreover, using Lemma \ref{lem:approx_sol}, we deduce that $W$ satisfies
\begin{eqnarray*}
\p_\tau W &=&\dv_x(x W) + \Delta_x w - R\dv_x\left((\alpha_1(z)-c)W\right)\\
&&- R^{N+1}\dv_x \left[ \tilde B_1 \left(z, \frac{U}{R^N}  \right) - \tilde B_1 \left(z, \frac{\Uapp}{R^N}  \right)\right]\\
&&+ U^\text{rem},
\end{eqnarray*}
with $R= e^\tau$, $z=R x + c \frac{R^2-1}{2}$,  and we recall that the remainder $U^{rem}$ satisfies
\be\label{in:Urem}
\|U^\text{rem}(\tau)\|_{L^\infty(\R^N)}+\|U^\text{rem}(\tau)\|_{L^2(e^{\gamma |x|^2})}\leq C e^{-\tau} 
\ee
for some $\gamma>0.$

Then, using the bounds on $U, \Uapp$ together with  the regularity assumptions on $\tilde B$, it can be easily proved that
$$
\tilde B_1 \left(z, \frac{U(\tau,x)}{R^N}  \right) - \tilde B_1 \left(z, \frac{\Uapp(\tau,x)}{R^N}  \right)= 2\alpha_2(z) f_0(z) \frac{F_M(x) W(\tau,x)}{R^{2N}} + b(\tau,x),
$$
and the function $b$ is such that there exists $C>0$ such that
$$
\forall (\tau,x)\in\R_+\times\R^N,\quad|b(\tau,x)|\leq C \left(\left|\frac{W(\tau,x)}{R^{N}}\right|^2 + R^{-2N-1} |W(\tau,x)|\right).
$$
We define a function $\tilde W$ by
$$
\tilde W(\tau, x)=W_0(x,z)+ e^{-\tau} W_1(x,z),
$$
with $W_0(x,z)=f_0(z) h_m(x)$ and 
$$
-\Delta_z W_1 + \dv_z (\alpha_1 W_1)= 2 \Delta_{yz} W_0 -\dv_x((\alpha_1-c) W_0) - \mathbf 1_{N=1}2\dv_z(\alpha_2 f_0 F_M W_0).
$$
Notice that by definition of $f_0$ and $c$, the compatibility condition is always satisfied, and
$$
W_1(x,z)= f_1(z)\cdot \nabla_y h_m (x)+ \mathbf 1_{N=1} w_1 (z) F_M (x) h_m(x),
$$
with
$$
-\Delta_z w_1 + \dv_z (\alpha_1 w_1)= -2 \dv_z(\alpha_2 f_0^2).
$$
Let $\tau_0>0$ such that 
$$
\tilde W(\tau,x)\geq \frac{1}{2}f_0(z) h_m(x)\quad\forall \tau\geq \tau_0,\ \forall y\in\R^N.
$$
For further purposes, we also choose $\tau_0$ such that
$$
\|W(\tau,\cdot)\|_1\leq 2 \ell(U_{ini})\quad\forall\tau\geq \tau_0.
$$
(Notice that if $\ell(U_{ini})=0$ there is nothing to prove).

Using calculations similar to the ones performed in the proof of Proposition \ref{prop:compactness}, we infer that for $\tau\geq \tau_0$,
\begin{eqnarray*}
&&\frac{d}{d\tau } \int_{\R^N} \left|\frac{W}{\tW}  \right|^2 \tW(\tau)+ \frac{m-N}{4}\int_{\R^N} \left|\frac{W}{\tW}  \right|^2 \tW(\tau) + 2 \int_{\R^N} \left|\nabla\frac{W}{\tW}  \right|^2 \tW(\tau)\\
&\leq & C \int_{\R^N} \left(\frac{W(\tau,x)}{\tW(\tau,x)}\right)^2 \frac{dx}{(1+|x|^2)^{1+\frac{m}{2}}} \\
&&+  C e^{(1-N)\tau}\int_{\R^N} |W(\tau,x)|^2 \left|  \nabla_y \frac{W(\tau,x)}{\tW(\tau,x)} \right|\:dx\\
&&+ \int_{\R^N} \left|\frac{W(\tau,x)}{\tW(\tau,x)}\right|\; \left|U^{rem}(\tau,x)  \right|\:dx.
\end{eqnarray*}
Using the same arguments as in the third step of the proof of Proposition \ref{prop:compactness}, we deduce that if $m>2(N+1)$,
\begin{eqnarray*}
 \int_{\R^N} \left(\frac{W(\tau,x)}{\tW(\tau,x)}\right)^2 \frac{dx}{(1+|x|^2)^{1+\frac{m}{2}}} &\leq &  \frac{m-N}{16}\int_{\R^N} \left|\frac{W(\tau)}{\tW(\tau)}  \right|^2 \tW(\tau) \\
&&+\frac{1}{2}\int_{\R^N} \left|\nabla\frac{W(\tau)}{\tW(\tau)}  \right|^2 \tW(\tau)\\
&&+ C \ell(U_{ini})^2.
\end{eqnarray*}
Similarly, the calculations of the fourth step in the proof of Proposition \ref{prop:compactness} yield
\begin{eqnarray*}
&&e^{(1-N)\tau}\int_{\R^N} |W(\tau,x)|^2 \left|  \nabla_x \frac{W(\tau,x)}{\tW(\tau,x)} \right|\:dx\\
&\leq & C \|W(\tau)\|_{L^1(\R^N)}^{1/N} \left[\int_{\R^N} \left|\frac{W(\tau)}{\tW(\tau)}  \right|^2 \tW(\tau)+\int_{\R^N} \left|\nabla\frac{W(\tau)}{\tW(\tau)}  \right|^2 \tW(\tau) \right]\\
&\leq & C \ell(U_{ini})^{1/N}  \left[\int_{\R^N} \left|\frac{W(\tau)}{\tW(\tau)}  \right|^2 \tW(\tau)+\int_{\R^N} \left|\nabla\frac{W(\tau)}{\tW(\tau)}  \right|^2 \tW(\tau) \right].
\end{eqnarray*}
Eventually, using the Cauchy-Schwarz inequality together with the bound \eqref{in:Urem}, we infer that
\begin{eqnarray*}
&& \int_{\R^N}\left|\frac{W(\tau,x)}{\tW(\tau,x)}\right|\; \left|U^{rem}(\tau,x)  \right|\:dx\\&\leq&\| U^{rem}(\tau)\|_{L^2(\tilde W(\tau)^{-1})}\left(\int_{\R^N} \left|\frac{W(\tau)}{\tW(\tau)}  \right|^2 \tW(\tau)\right)^{1/2}\\
&\leq & C\| U^{rem}(\tau)\|_{L^2(e^{\gamma|x|^2})}\left(\int_{\R^N} \left|\frac{W(\tau)}{\tW(\tau)}  \right|^2 \tW(\tau)\right)^{1/2}\\
&\leq & C e^{-\tau}\left(\int_{\R^N} \left|\frac{W(\tau)}{\tW(\tau)}  \right|^2 \tW(\tau)\right)^{1/2}\\
&\leq & C + \frac{m-N}{16}\int_{\R^N} \left|\frac{W(\tau)}{\tW(\tau)}  \right|^2 \tW(\tau).
\end{eqnarray*}

Gathering all the terms, we deduce that there exists a constant $C_m$, depending only on $N$ and $m$, such that if $\ell(U_{ini})\leq C_m$, then for all $\tau\geq \tau_0,$
$$
\frac{d}{d\tau } \int_{\R^N} \left|\frac{W(\tau)}{\tW(\tau)}  \right|^2 \tW(\tau)+ \frac{m-N}{16}\int_{\R^N} \left|\frac{W(\tau)}{\tW(\tau)}  \right|^2 \tW(\tau) +  \int_{\R^N} \left|\nabla\frac{W(\tau)}{\tW(\tau)}  \right|^2 \tW(\tau)\leq C .
$$
Compactness of a subsequence $W(\tau_n)$ follows. Hence the $\omega$-limit set is non-empty, and thus $\ell(U_{ini})=0$.

\section*{Appendix A}

\noindent \textbf{Lemma A.1.}\textit{ Assume that the flux $A$ satisfies \eqref{hyp:A1}, \eqref{hyp:A2}.
Let $v\in W^{1,\infty}(\T^N)$ be a periodic stationary solution of \eqref{LCS}, and let $u\in L^\infty_\text{loc}([0,\infty), L^\infty(\R^N))\cap \mathcal C([0,\infty), L^1_\text{loc}(\R^N))$ be the unique solution of \eqref{LCS} with initial data $u_{ini}\in v(y) + L^1\cap L^\infty(\R^N).$ Then $u \in L^\infty([0,\infty)\times \R^N).$}

\begin{proof}
This result was proved in \cite{longtime} in the case $N=1$. When $N\geq 2$, the proof goes along the same lines; the only difference lies in the use of the Sobolev embeddings, which depend on the dimension.
Hence we merely recall here the main steps of the proof, with an emphasis on  the case $N\geq 2$.

In the rest of the proof, we set $f(t,y)=u(t,y)-v(y).$ Then $f$ solves the equation
\be\label{eq:f}
\p_t f + \dv_y B(y,f) - \Delta_y f=0,
\ee
and according to \eqref{hyp:A2} the flux $B$ is such that for all $f\in\R$, 
$$
\begin{aligned}
\left| \dv_y B(y,f) \right|\leq C (|f| + |f|^n),\\
\left| \p_f  B(y,f) \right|\leq C (|f| + |f|^n), \end{aligned}
$$
where the exponent $n$ is such that $n<(N+2)/N$.
Moreover, 
$$
\|f(t)\|_{L^1(\R^N)}\leq \|u_{ini} - v\|_{L^1(\R^N)}\quad\forall t\geq 0.
$$
For $q\geq 1$ arbitrary, multiply \eqref{eq:f} by $|f|^q$, and integrate over $\R^N.$ Using a few integrations by parts (see \cite{longtime}), we are led to 
\be\label{in:f}
\frac{d}{dt}\int_{\R^N}|f|^{q+1} + c_q \int_{\R^N}\left|\nabla_y |f|^\frac{q+1}{2}  \right|^2\leq C_q\left(\int_{\R^N} |f|^{q+1} + \int_{\R^N} |f|^{q+n}\right)
\ee
We then use Sobolev embeddings in order to control the $L^{q+1}$ and $L^{q+n}$ norms in the right-hand side. We distinguish between the cases $N=2$ and $N\geq 3$, since the space $H^1$ is critical in dimension two.

$\bullet$ If $N=2$, then $H^1(\R^2)\subset L^p(\R^2)$ for all $ p\in[2,\infty)$. Interpolating $L^{q+n}$ between $L^1$ and $L^p$ for some $p$ sufficiently large, we have
\begin{eqnarray*}
\|f\|_{L^{q+n}(\R^2)}&\leq & \|f\|_1^\theta \|f\|_p^{1-\theta}\quad \text{with } \frac{1}{q+n}= \frac{\theta}1 + \frac{1-\theta} p\\
&\leq & \|f\|_1^\theta \left\| |f|^\frac{q+1}{2} \right\|_{\frac{2p}{q+1}}^\frac{2(1-\theta)}{q+1}\\
&\leq &C_p \|f\|_1^\theta \left\| \nabla |f|^\frac{q+1}{2} \right\|_2^\frac{2(1-\theta)}{q+1}.
\end{eqnarray*}
Notice that
$$
\frac{q+n}{q+1}(1-\theta)=\frac{q+n-1}{(q+1)\left(1-\frac{1}p\right)},
$$
and 
$$
\frac{q+n-1}{q+1}<1\quad \forall q\geq 1
$$
since $n<(N+2)/N$. Thus, we choose $p>1$ such that 
$$
\frac{q+n-1}{(q+1)\left(1-\frac{1}p\right)}<1.
$$
Young's inequality then implies that for all $\lambda>0,$ there exists a constant $C_{\lambda,q}$ and an exponent $q_1$ such that
\be\label{in:q+n}
\int_{\R^N}|f|^{q+n}\leq \lambda \left\| \nabla |f|^\frac{q+1}{2} \right\|_2^{2} + C_{\lambda,q}\| f\|_1^{q_1} .
\ee
The other term in the right-hand side of \eqref{in:f} can be bounded in a similar fashion: we have, for all $\lambda>0$,
\be\label{in:Young}
\int_{\R^N}|f|^{q+1}\leq \lambda \left\| \nabla |f|^\frac{q+1}{2} \right\|_2^{2} + C_{\lambda,q}\| f\|_1^{q_2},
\ee
for some exponent $q_2$ which can be explicitely computed.
 Choosing an appropriate parameter $\lambda$, we infer that there exist $q_1,q_2>0$ such that
$$
\frac{d}{dt}\int_{\R^N}|f|^{q+1} + c_q \int_{\R^N}\left|\nabla_y |f|^\frac{q+1}{2}  \right|^2\leq C_q(\|f\|_1^{q_1} +\|f\|_1^{q_2}  ).
$$
Using \eqref{in:Young} one more time leads to 
$$
\frac{d}{dt}\int_{\R^N}|f|^{q+1} + c_q \int_{\R^N}|f|^{q+1} \leq C_q(\|f\|_1^{q_1} +\|f\|_1^{q_2}  )\leq C.
$$
Using a Gronwall-type argument, we infer that $f\in L^\infty([0,\infty), L^\frac{q+1}{2}(\R^N))$ for all $q\geq 1$.
\vskip1mm

$\bullet$ When $N\geq 3$, we use the Sobolev embedding $H^1(\R^N)\subset L^{p^*}(\R^N)$, where 
$$
p^*=\frac{2N}{N-2}.
$$
Interpolating $L^{q+n}$ between $L^1$ and $L^\frac{p^*(q+1)}{2}$, we obtain
\begin{eqnarray*}
\|f\|_{q+n}\leq \|f\|_\frac{p^*(q+1)}{2}^\theta \|f\|_1^{1-\theta}
&\leq&\left\| |f|^\frac{q+1}{2}\right\|_{p^*}^\frac{2\theta}{q+1} \|f\|_1^{1-\theta}\\
&\leq & C \left\|\nabla |f|^\frac{q+1}{2}\right\|_2^\frac{2\theta}{q+1} \|f\|_1^{1-\theta},
\end{eqnarray*}
where the parameter $\theta\in(0,1)$ is given by
$$
\frac{1}{q+n}=\frac{2\theta}{p^* (q+1)} + \frac{1-\theta}{1}.
$$ 
It can be checked that
$$
n<\frac{N+2}{N}\Rightarrow \frac{\theta(q+n)}{q+1}<1.
$$
Hence \eqref{in:q+n} holds when $N\geq 3$. Inequality \eqref{in:Young} is proved with similar arguments. As in the two-dimensional case, we deduce that $f\in L^\infty([0,\infty), L^q(\R^N))$ for all $q$. Using Theorem 8.1 in Chapter III of \cite{ladypara} (see \cite{longtime} for details), we infer eventually that $f\in L^\infty([0,\infty)\times\R^N)$.

\end{proof}

\section*{Appendix B}

\noindent \textbf{Lemma A.2.}\textit{Let $M>M'$ be arbitrary. Then
$$
F_M(y)> F_{M'}(y)\quad \forall y\in\R^N.
$$
As a consequence,
$$
\| F_M - F_{M'}\|_1=M-M'.
$$
}

\begin{proof}
 The arguments are exactly the ones which lead to the uniqueness of stationary solutions of \eqref{eq:h_Ngrand}, \eqref{eq:h_N=1}, and they can be found in \cite{AEZ}. We recall the main steps below for the reader's convenience. 

Let $F:= F_M-F_{M'}$. Then $F\in L^1\cap \mathcal C^2(\R^N)$, and $\int_{\R^N} F>0$. Hence the set 
$$
\Theta:=\{x\in\R^N, F(x)>0\}
$$
is non-empty. The idea is to prove that $F_+=F\mathbf 1_\Theta$ satisfies a linear elliptic equation; since $F_+\geq 0$, $F_+$ cannot vanish anywhere, and thus $F_+(x)>0$ for all $x\in\R^N$.

Let us now derive an equation on $F_+$. Substracting the equations on $F_M$ and $F_{M'}$, we have
$$
-\sum_{1\leq i,j\leq N}\eta_{i,j}\frac{\p^2 F}{\p x_i \p x_j} + \dv_x(bF)=0,
$$
where
$$
b(x)=a(F_M(x) + F_{M'}(x)) - x,\quad x\in\R^N;
$$
notice that $a=0$ if $N\geq 2$. Since $F\in H^2(\R^N)$, we have 
$$
\dv_x(bF)\mathbf 1_\Theta= \dv_x(bF_+)
$$
almost everywhere. Thus, we obtain
$$
-\sum_{1\leq i,j\leq N}\eta_{i,j}\mathbf 1_\Theta\frac{\p^2 F}{\p x_i \p x_j} + \dv_x(bF_+)=0.
$$
Integrating the above equation on $\R^N$ leads to 
$$
\int_{\Theta} \sum_{1\leq i,j\leq N}\eta_{i,j}\frac{\p^2 F}{\p x_i \p x_j}=0.
$$
Let us now perform the change of variables \eqref{CV}, which changes the  matrix $\eta$ into identity: setting $\tilde F(y)=F(Py)$, and $\tilde \Theta:=\{\tilde F>0\}$, we infer
$$
\int_{\tilde \Theta}\Delta_y \tilde F= C\int_{\Theta} \sum_{1\leq i,j\leq N}\eta_{i,j}\frac{\p^2 F}{\p x_i \p x_j}=0.
$$
Moreover, $\tilde F\in H^2\cap W^{2,1}(\R^N)$, and thus Lemma 7 in \cite{AEZ} applies. We deduce that 
$$
\Delta_y(\tilde F \mathbf 1_{\tilde \Theta})= \mathbf1_{\tilde \Theta} \Delta_y \tilde F,
$$
and thus 
$$
\sum_{1\leq i,j\leq N}\eta_{i,j}\mathbf 1_\Theta\frac{\p^2 F}{\p x_i \p x_j}=\sum_{1\leq i,j\leq N}\eta_{i,j}\frac{\p^2 F_+}{\p x_i \p x_j}.
$$
Eventually, $F_+$ solves the elliptic equation
$$
-\sum_{1\leq i,j\leq N}\eta_{i,j}\frac{\p^2 F_+}{\p x_i \p x_j}+\dv_x(bF_+)=0,
$$
with $b\in L^\infty_\text{loc}(\R^N)$. Using either a unique continuation principle or Harnack's inequality (see \cite{GT}, Theorem 8.20), we infer that if $F_+$ vanishes at some point $x$ in $\R^N$, then $F_+$ is identically zero on $\R^N$, which is absurd. Hence $F_+(x)>0$ for all $x\in\R^N$, and thus $\R^N\setminus \Theta=\emptyset,$ which means that $F(x)>0$ for all $x\in\R^N$.

\end{proof}

\section*{Acknowledgements}

I wish to thank Adrien Blanchet, Jean Dolbeault, and Michal Kowalczyk, for very fruitful and stimulating discussions.

\bibliography{ratchets}

\end{document}